# Some observations about determinants which are connected with Catalan numbers and related topics

*Johann Cigler*

Fakultät für Mathematik, Universität Wien

**Abstract**

In this (mostly expository) paper I want to share some observations prompted by a class of matrices whose determinants are Catalan numbers. Considering different methods of proof we obtain some generalizations and $q-$analogues and connections with Hankel determinants. Finally we state some conjectures.

## 1. Introduction

This (mostly expository) paper has an unusual genesis. In an internet posting an amateur mathematician, Tony Foster, looking for interesting properties of Pascal's triangle, observed that the elements of the "diagonals" $\{1,1\}, \{1,3,1\}, \{1,6,5,1\}, \cdots,$ of Pascal's triangle whose sums are Fibonacci numbers $1+1=2=F_3,\ 1+3+1=5=F_5,$ $1+6+5+1=13=F_7, \cdots$, yield matrices whose determinants are Catalan numbers:

$$\det(1) = 1 = C_1,\quad \det\begin{pmatrix} 1 & 1 \\ 1 & 3 \end{pmatrix} = 2 = C_2,\quad \det\begin{pmatrix} 1 & 1 & 0 \\ 1 & 3 & 1 \\ 1 & 6 & 5 \end{pmatrix} = 5 = C_3,$$

$$\det\begin{pmatrix} 1 & 1 & 0 & 0 \\ 1 & 3 & 1 & 0 \\ 1 & 6 & 5 & 1 \\ 1 & 10 & 15 & 7 \end{pmatrix} = 14 = C_4 \quad \text{and more generally}$$

$$\det\left(\binom{i+j+1}{2j}\right)_{i,j=0}^{n-1} = \det\left(\binom{i+j+1}{i-j+1}\right)_{i,j=0}^{n-1} = C_n = \frac{1}{n+1}\binom{2n}{n} \tag{1}$$

for non-negative integers $n$.

For example elementary column operations give

$$\det\begin{pmatrix} 1 & 1 & 0 & 0 \\ 1 & 3 & 1 & 0 \\ 1 & 6 & 5 & 1 \\ 1 & 10 & 15 & 7 \end{pmatrix} = \det\begin{pmatrix} 1 & 1-1 & 0 & 0 \\ 1 & 3-1 & 1 & 0 \\ 1 & 6-1 & 5 & 1 \\ 1 & 10-1 & 15 & 7 \end{pmatrix} = \det\begin{pmatrix} 1 & 0 & 0 & -\frac{0}{2} & 0 \\ 1 & 2 & 1 & -\frac{2}{2} & 0 \\ 1 & 5 & 5 & -\frac{5}{2} & 1 \\ 1 & 9 & 15 & -\frac{9}{2} & 7 \end{pmatrix}$$

$$= \det\begin{pmatrix} 1 & 0 & 0 & 0 \\ 1 & 2 & 0 & 0 \\ 1 & 5 & \frac{5}{2} & 1-\frac{2}{5}\frac{5}{2} \\ 1 & 9 & \frac{21}{2} & 7-\frac{2}{5}\frac{21}{2} \end{pmatrix} = \det\begin{pmatrix} 1 & 0 & 0 & 0 \\ 1 & 2 & 0 & 0 \\ 1 & 5 & \frac{5}{2} & 0 \\ 1 & 9 & \frac{21}{2} & \frac{14}{5} \end{pmatrix} = 1 \cdot 2 \cdot \frac{5}{2} \cdot \frac{14}{5} = \frac{C_1}{C_0} \cdot \frac{C_2}{C_1} \cdot \frac{C_3}{C_2} \cdot \frac{C_4}{C_3} = C_4.$$

This caught my attention. I sought proofs which in turn led to some other interesting results. Probably some of them are known, but I could not find them in the literature. Relevant references are very welcome. In order to make the paper readable by non-specialists I do not hesitate to recall proofs of well-known results.

It turns out that the identities (1) are equivalent with the following identities for Catalan numbers

$C_0 = 1$, $C_1 - C_0 = 0$, $C_2 - 3C_1 + C_0 = 0$, $C_3 - 5C_2 + 6C_1 - C_0 = 0$,
$C_4 - 7C_3 + 15C_2 - 10C_1 + C_0 = 0$, and more generally with

$$\sum_{j=0}^{n} (-1)^{n-j} \binom{n+j}{n-j} C_j = [n=0]. \tag{2}$$

Thus at first I looked for proofs of these identities.

**1.1.** The connection with Fibonacci numbers suggests the following interpretation of the left-hand side of (2): Let $F_n(x) = \sum_{j=0}^{\lfloor n/2 \rfloor} (-1)^j \binom{n-j}{j} x^{n-2j}$ be Fibonacci polynomials and let $\Lambda$ be the linear functional on the polynomials defined by $\Lambda(F_n(x)) = [n=0]$. Then the Fibonacci polynomials are orthogonal with respect to $\Lambda$ with $\Lambda(x^{2n}) = C_n$. Since

$$F_{2n}(x) = \sum_{j=0}^{n} (-1)^j \binom{2n-j}{j} x^{2n-2j} = \sum_{j=0}^{n} (-1)^{n-j} \binom{n+j}{n-j} x^{2j}$$ we see that (2) is equivalent with

$\Lambda(F_{2n}(x)) = [n=0]$.

**1.2.** Other methods which lead to (2) are expansions of the polynomials $\binom{x-1}{n}$ with respect to Gould polynomials $G_n(x) = \frac{x}{2n+x}\binom{2n+x}{n}$ or hypergeometric summation formulae. A combinatorial proof of a more general formula has been given by Ricky X. F. Chen [3], Theorem 1.1.

These methods also lead to

$$\det\left(\binom{i+j+k}{i-j+1}\right)_{i,j=0}^{n-1} = \det\left(\binom{j+k}{i-j+1}\right)_{i,j=0}^{n-1} = C_n^{(k)} \qquad (3)$$

where the numbers

$$C_n^{(k)} = \frac{k}{n+k}\binom{2n+k-1}{n} = \frac{k}{2n+k}\binom{2n+k}{n} = \binom{2n+k-2}{n} - \binom{2n+k-2}{n-2} \qquad (4)$$

are the coefficients of the $k-$th power of the generating function of the Catalan numbers,

$$\sum_{n\geq 0} C_n^{(k)} z^n = \left(\sum_{n\geq 0} C_n z^n\right)^k.$$

**1.3** A different method to prove (1) uses a result by Christian Krattenthaler about determinants of binomial coefficients and establishes a surprising relation with Hankel determinants:

$$\det\left(\binom{i+j+m}{i-j+m}\right)_{i,j=0}^{n-1} = \det\left(C_{n+i+j}\right)_{i,j=0}^{m-1}. \qquad (5)$$

The appearance of the Hankel determinants $\det\left(C_{n+i+j}\right)_{i,j=0}^{m-1}$ at first glance seemed to be a lucky coincidence, but some experimentation led to the conjecture that for arbitrary monic orthogonal polynomials $p_n(x) = \sum_{j=0}^{n}(-1)^{n-j} p(n,j) x^j$ with moments $M_n$ the identity

$$\det\left(p(i+m,j)\right)_{i,j=0}^{n-1} = \frac{\det\left(M_{n+i+j}\right)_{i,j=0}^{m-1}}{\det\left(M_{i+j}\right)_{i,j=0}^{m-1}} \qquad (6)$$

holds. In the mean-time this conjecture has been proved by M. Tyson [18].

**1.4.** It is also natural to look for $q-$analogues. As usual we write $[n] = [n]_q = \dfrac{1-q^n}{1-q}$,

$[n]! = [1][2]\cdots[n]$ and $\begin{bmatrix} n \\ k \end{bmatrix} = \begin{bmatrix} n \\ k \end{bmatrix}_q = \dfrac{[n]!}{[k]![n-k]!}$ for $0 \leq k \leq n$. Further we let

$(x;q)_n = (1-x)(1-qx)\cdots(1-q^{n-1}x)$.

As analogues of (3) we get

$$\det\left(q^{\binom{i-j}{2}}\begin{bmatrix} i+j+k \\ i-j+1 \end{bmatrix}\right)_{i,j=0}^{n-1} = C_n^{(k)}(q) \tag{7}$$

and a more intricate identity

$$\det\left(q^{2\binom{i-j+1}{2}}\frac{\left(-q^{j+k};q\right)_{i-j+1}}{\left(-q;q\right)_{i-j+1}}\begin{bmatrix} j+k \\ i-j+1 \end{bmatrix}\right)_{i,j=0}^{n-1} = \frac{\left(-q^{n+1};q\right)_{k-1}}{\left(-q;q\right)_{k-1}} C_n^{(k)}(q) \tag{8}$$

where $C_n^{(k)}(q) = \dfrac{[k]}{[2n+k]}\begin{bmatrix} 2n+k \\ n \end{bmatrix} = \begin{bmatrix} 2n+k-2 \\ n \end{bmatrix} - q^k \begin{bmatrix} 2n+k-2 \\ n-2 \end{bmatrix}$ is a $q-$analogue of $C_n^{(k)}$.

**1.5.** In another direction we get

$$\det\left(\binom{i+j+k+m}{i-j+m}\right)_{i,j=0}^{n-1} = \det\left(\binom{j+k+m}{i-j+m}\right)_{i,j=0}^{n-1} = \det\left(C_{n-i+j}^{(2i+k+1)}\right)_{i,j=0}^{m-1} \tag{9}$$

and

$$\det\left(\frac{x+2i-1+2m}{x+i+j+m-1}\binom{x+i+j+m-1}{i-j+m}\right)_{i,j=0}^{n-1} = \det\left(\binom{2n+2j+x-1}{n-i+j}\right)_{i,j=0}^{m-1} \tag{10}$$

and their $q-$ analogues.

Of course the above results remain true if we reduce them modulo a prime $p$. But what happens when we consider the residues modulo $p$ not as elements of the residue ring but as integers or real numbers? There seems to be no reason to expect nice results. But in some special cases nevertheless nice results occur.

It turns out that

$$\det\left(\binom{i+j+1}{i-j+1} \bmod 2\right)_{i,j=0}^{n-1} = C_n \bmod 2 \tag{11}$$

remains true with this interpretation. And computations suggest that

$$\det\left(\binom{i+j+k}{i-j+1} \bmod 2\right)_{i,j=0}^{n-1} = (-1)^{k-1} C_n^{(k)} \bmod 2, \tag{12}$$

and

$$\det\left(\binom{i+j+m}{i-j+m} \bmod 2\right)_{i,j=0}^{n-1} = \det\left(C_{n-i+j}^{(2i+1)} \bmod 2\right)_{i,j=0}^{m-1}. \tag{13}$$

Another nice result seems to occur for $p = 3$. Here we get

$$\det(1) = 1, \quad \det\begin{pmatrix} 1 & 1 \\ 1 & 0 \end{pmatrix} = -1, \quad \det\begin{pmatrix} 1 & 1 & 0 \\ 1 & 0 & 1 \\ 1 & 0 & 2 \end{pmatrix} = -1, \quad \det\begin{pmatrix} 1 & 1 & 0 & 0 \\ 1 & 0 & 1 & 0 \\ 1 & 0 & 2 & 1 \\ 1 & 1 & 0 & 1 \end{pmatrix} = -1,$$

and more generally we guess that

$$\det\left(\binom{i+j+1}{i-j+1} \bmod 3\right)_{i,j=0}^{n-1} = \mu(C_n), \tag{14}$$

where $n \bmod 3 \in \{0,1,2\}$ and $\mu(n)$ denotes the residue modulo $3$ where the residues belong to $\{-1,0,1\}$, i.e. $\mu(3n) = 0$, $\mu(3n+1) = 1$ and $\mu(3n+2) = -1$.

I want to thank Christian Krattenthaler, Michael Schlosser and Darij Grinberg for helpful informations.

## 2. Different proofs

The matrices $\left(\binom{i+j+1}{i-j+1}\right)_{i,j=0}^{n-1}$ are "almost triangular" in the sense that all entries $a(i,j)$ with $j > i+1$ vanish. Our first proof uses elementary matrix operations.

**Lemma 1**

*Let $T = \big(t(i,j)\big)_{i,j\geq 0}$ be a triangular matrix with $t(i,i) = 1$ for all $i$ and let $T_{n,1} = \big(t(i+1,j)\big)_{i,j=0}^{n-1}$. If there are numbers $M_n \neq 0$ such that*

$$\sum_{j=0}^{n}(-1)^{n-j} t(n,j) M_j = [n=0] \tag{15}$$

*then*

$$\det\big(T_{n,1}\big) = \det\big(t(i+1,j)\big)_{i,j=0}^{n-1} = M_n. \tag{16}$$

**Proof**

Let $t_j$ be column $j$ of $T_{n,1} = \bigl(t(i+1,j)\bigr)_{i,j=0}^{n-1}$ and let $S_n$ be the matrix with columns

$$s_j = \bigl(s(i,j)\bigr)_{i=0}^{n-1} = t_j + \frac{1}{M_j}\sum_{k=1}^{j}(-1)^k M_{j-k} t_{j-k}.$$

Then we get $M_j s_j = \sum_{k=0}^{j}(-1)^{j-k} M_k t_k = M_j t_j - M_{j-1} s_{j-1}$.

For each $i$ this implies

$$s(i,j)M_j = \sum_{k=0}^{j}(-1)^k M_{j-k} t(i+1, j-k) = 0$$

for $j = i+1$ by (15) and $s(i,j)M_j = M_j t(i,j) - M_{j-1} s(i, j-1) = 0$ for $j > i+1$ by induction. Thus $S_n$ is triangular.

Furthermore we get

$$s(i,i)M_i = -\sum_{k=0}^{i+1}(-1)^k M_{i+1-k} t(i+1, i+1-k) + M_{i+1} = M_{i+1}.$$

Thus $s(i,i) = \dfrac{M_{i+1}}{M_i}$.

Since $S_n$ is derived from $T_{n,1}$ by elementary column operations we have

$$\det T_{n,1} = \det S_n = \prod_{i=0}^{n-1}\frac{M_{i+1}}{M_i} = M_n. \qquad \square$$

This simple observation has many applications.

For example for $t(i,j) = \binom{i}{j}$ (15) reduces to $\sum_{k=0}^{n}(-1)^{n-k}\binom{n}{k}\cdot 1 = \bigl(1-1\bigr)^n = [n=0]$

which implies the well-known result $\det\left(\binom{i+1}{j}\right)_{i,j=0}^{n-1} = 1$.

For the matrices (1) we choose

$$T = A = \big(a(i,j)\big)_{i,j\geq 0} = \left(\binom{i+j}{i-j}\right)_{i,j\geq 0} = \begin{pmatrix} 1 & 0 & 0 & 0 & 0 & \cdots \\ 1 & 1 & 0 & 0 & 0 & \cdots \\ 1 & 3 & 1 & 0 & 0 & \cdots \\ 1 & 6 & 5 & 1 & 0 & \cdots \\ 1 & 10 & 15 & 7 & 1 & \cdots \\ \vdots & \vdots & \vdots & \vdots & \vdots & \ddots \end{pmatrix},$$

In this case identity (15) reduces to the identity (2) which implies (1).

**Remark**

Let $\big(b(i,j)\big)_{i,j\geq 0}$ be the inverse matrix of $\left(\binom{i+j}{i-j}\right)_{i,j\geq 0}$. Then (2) is equivalent with

$b(n,0) = (-1)^n C_n$. Note that by Cramer's rule $b(n,0) = (-1)^n \det\left(\binom{i+j+1}{i-j+1}\right)_{i,j=0}^{n-1}$.

### 2.1. Some proofs of the identities for Catalan numbers

#### 2.1.1. Orthogonal polynomials

As already mentioned identity (2) occurs in a natural way in the computation of the Hankel determinants $\det\big(C_{i+j}\big)_{i,j=0}^{n-1}$ of the Catalan numbers with the method of orthogonal polynomials.

Let me recall some relevant facts.

**Lemma 2.**

*Let $(a(n))$ be a sequence of real numbers with $a(0)=1$.*

*If all Hankel determinants $\det(a(i+j))_{i,j=0}^{n-1} \neq 0$ then the polynomials*

$$p_n(x) = \frac{1}{\det(a(i+j))_{i,j=0}^{n-1}} \det\begin{pmatrix} a(0) & a(1) & \cdots & a(n-1) & 1 \\ a(1) & a(2) & \cdots & a(n) & x \\ a(2) & a(3) & \cdots & a(n+1) & x^2 \\ \vdots & & & & \vdots \\ a(n) & a(n+1) & \cdots & a(2n-1) & x^n \end{pmatrix} \quad (17)$$

*are orthogonal with respect to the linear functional $\Lambda$ defined by*

$$\Lambda(x^n) = a(n). \qquad (18)$$

*This means that* $\Lambda\left(p_n(x)p_m(x)\right) = 0$ *if* $m \neq n$ *and* $\Lambda\left(p_n^2(x)\right) \neq 0$.

*In particular for* $m = 0$ *we get*

$$\Lambda\left(p_n(x)\right) = [n = 0]. \tag{19}$$

*By Favard's theorem there exist numbers* $s(n), t(n)$ *such that*

$$p_n(x) = (x - s(n-1))p_{n-1}(x) - t(n-2)p_{n-2}(x). \tag{20}$$

*If on the other hand for given* $s(n)$ *and* $t(n)$ *we define* $c(n, j)$ *by*

$$\begin{aligned} c(0, j) &= [j = 0] \\ c(n, 0) &= s(0)c(n-1, 0) + t(0)c(n-1, 1) \\ c(n, j) &= c(n-1, j-1) + s(j)c(n-1, j) + t(j)c(n-1, j+1) \end{aligned} \tag{21}$$

*and* $p_n(x)$ *by (20) and* $\Lambda$ *by (19) then* $c(n, 0) = \Lambda\left(x^n\right)$ *and*

$$\sum_{k=0}^{n} c(n, k) p_k(x) = x^n. \tag{22}$$

*This implies*

$$c(n, k) = \frac{\Lambda\left(x^n p_k(x)\right)}{\Lambda\left(p_k(x)^2\right)}. \tag{23}$$

Let $p_n(x) = \sum_{j=0}^{n} (-1)^{n-j} p(n, j) x^j$.

From $\dfrac{\Lambda\left(p_k(x)p_{n+k}(x)\right)}{\Lambda\left(p_k(x)^2\right)} = [n = 0]$ we get

$$\sum_{j=0}^{n} (-1)^{n-j} p(n+k, j+k) c(j+k, k) = [n = 0] \tag{24}$$

because

$$\begin{aligned}
\frac{\Lambda\left(p_k(x)p_{n+k}(x)\right)}{\Lambda\left(p_k(x)^2\right)} &= \sum_{j=0}^{n+k} (-1)^{n+k-j} p(n+k, j) \frac{\Lambda\left(x^j p_k(x)\right)}{\Lambda\left(p_k(x)^2\right)} \\
&= \sum_{j=k}^{n+k} (-1)^{n+k-j} p(n+k, j) \frac{\Lambda\left(x^j p_k(x)\right)}{\Lambda\left(p_k(x)^2\right)} = \sum_{j=k}^{n+k} (-1)^{n+k-j} p(n+k, j) c(j, k) \\
&= \sum_{j=0}^{n} (-1)^{n-j} p(n+k, j+k) c(j+k, k).
\end{aligned}$$

Therefore we get (24) which by (16) implies

$$\det\left(p(i+k+1, j+k)\right)_{i,j=0}^{n-1} = c(n+k, k). \tag{25}$$

As a simple example let us consider the polynomials

$$p_n(x) = \sum_{j=0}^{n} (-1)^{n-j} \begin{bmatrix} n \\ j \end{bmatrix} q^{(n-1)(n-j)} x^j \tag{26}$$

which are a $q-$analogue of $(x-1)^n$. They are orthogonal with respect to the linear functional $\Lambda_q$ defined by $\Lambda_q\left(x^n\right) = q^{\binom{n}{2}}$.

The orthogonality follows from the well-known formula

$$(x;q)_n = (1-x)(1-qx)\cdots(1-q^{n-1}x) = \sum_{j=0}^{n} (-1)^j q^{\binom{j}{2}} \begin{bmatrix} n \\ j \end{bmatrix} x^j,$$

which implies

$$\Lambda_q\left(x^m p_n(x)\right) = \sum_{j=0}^{n} (-1)^{n-j} \begin{bmatrix} n \\ j \end{bmatrix} q^{(n-1)(n-j)} q^{\binom{j+m}{2}} = (-1)^n q^{\binom{m}{2}+n^2-n} \sum_{j=0}^{n} (-1)^j \begin{bmatrix} n \\ j \end{bmatrix} q^{\binom{j}{2}} \left(\frac{1}{q^{n-m-1}}\right)^j$$

$$= \left(1 - \frac{1}{q^{n-m-1}}\right)\left(1 - \frac{1}{q^{n-m-2}}\right)\cdots\left(1 - \frac{1}{q^{-m}}\right) = 0$$

for $m < n$ and

$$\Lambda_q\left(x^n p_n(x)\right) = (1-q)(1-q^2)\cdots(1-q^n) \neq 0.$$

Note that for $q=1$ the polynomials $(x-1)^n$ are not orthogonal with respect to the linear functional $\Lambda_1$ since $\Lambda_1\left((x-1)^n(x-1)^n\right) = 0$.

The corresponding numbers $c(n,k) = \dfrac{\Lambda_q\left(x^n p_k(x)\right)}{\Lambda_q\left(x^k p_k(x)\right)}$ are

$$c(n,k) = q^{\binom{n}{2}-\binom{k}{2}} \frac{(1-q^{n+1-k})(1-q^{n+2-k})\cdots(1-q^n)}{(1-q)(1-q^2)\cdots(1-q^k)} = q^{\binom{n}{2}-\binom{k}{2}} \begin{bmatrix} n \\ k \end{bmatrix}.$$

From (25) we get

$$\det\left(q^{(i+k)(i+1-j)} \begin{bmatrix} i+k+1 \\ j+k \end{bmatrix}\right)_{i,j=0}^{n-1} = q^{\binom{n+k}{2}-\binom{k}{2}} \begin{bmatrix} n+k \\ k \end{bmatrix},$$

which can be simplified to

$$\det\left(q^{\binom{i-j}{2}}\begin{bmatrix}i+1+k\\j+k\end{bmatrix}\right)_{i,j=0}^{n-1}=\begin{bmatrix}n+k\\k\end{bmatrix}, \tag{27}$$

which is a well-known result.

Consider now the special Fibonacci polynomials

$$F_n(x)=\sum_{j=0}^{\lfloor\frac{n}{2}\rfloor}(-1)^j\binom{n-j}{j}x^{n-2j} \tag{28}$$

which satisfy $F_n(x)=xF_{n-1}(x)-F_{n-2}(x)$ with initial values $F_0(x)=1$ and $F_1(x)=x$.

Here we have $s(n)=0$ and $t(n)=1$. It is well known that the linear functional $\Lambda$ defined by $\Lambda(F_n(x))=0$ for $n>0$ and $\Lambda(F_0(x))=\Lambda(1)=1$ has as moments $\Lambda(x^{2n})=c(2n,0)=C_n$ and $\Lambda(x^{2n+1})=0$. Moreover we have (cf. e.g. [6])

$$c(2n+k,k)=\frac{k+1}{2n+k+1}\binom{2n+k+1}{n}=C_n^{(k+1)}. \tag{29}$$

Let us give a direct proof. The generating function $c(z)=\sum_{n\geq 0}C_n z^n$ of the Catalan numbers satisfies $c(z)=1+zc(z)^2$ which gives $c(z)^k=\sum_{n\geq 0}C_n^{(k)}z^n=c(z)^{k-1}+z^2c(z)^{k+1}$ and thus

$$C_n^{(k)}=C_n^{(k-1)}+C_{n-1}^{(k+1)}. \tag{30}$$

We need only compare this with (21) to get (29).

Applying $\Lambda$ to $F_{2n}(x)=\sum_{j=0}^{n}(-1)^j\binom{2n-j}{j}x^{2n-2j}=\sum_{j=0}^{n}(-1)^{n-j}\binom{n+j}{n-j}x^{2j}$ we get (2).

Since $F_{2n+k}(x)=\sum_{0\leq 2j+k\leq 2n+k}(-1)^{n-j}\binom{n+k+j}{n-j}x^{2j+k}$

we get $\dfrac{\Lambda(F_{2n+k}(x)F_k(x))}{\Lambda(F_k(x)^2)}=\sum_{-0\leq 2j+k\leq 2n+k}(-1)^{n-j}\binom{n+k+j}{n-j}\dfrac{\Lambda(x^{2j+k}F_k(x))}{\Lambda(F_k(x)^2)}.$

Therefore we get by (23) and $\Lambda(x^{2j+k}F_k(x))=0$ for $0\leq 2j+k<k$

$$\sum_{j=0}^{n}(-1)^{n-j}\binom{n+k+j}{n-j}C_j^{(k+1)} = [n=0] \tag{31}$$

which implies

$$\det\left(\binom{i+j+k}{i-j+1}\right)_{i,j=0}^{n-1} = C_n^{(k)}. \tag{32}$$

From

$$x^{n-1}F_{n+k+1}(x) = \sum_{0 \leq 2j \leq n+k+1}(-1)^j \binom{n+k+1-j}{j} x^{2n+k-2j}$$

$$= \sum_j (-1)^{n-j}\binom{k+1+j}{n-j} x^{k+2j}$$

we get

$$\sum_j (-1)^{n-j}\binom{k+1+j}{n-j}C_j^{(k+1)} = \frac{\Lambda\left(x^{n-1}F_{n+k+1}(x)F_k(x)\right)}{\Lambda\left(F_k(x)^2\right)} = [n=0] \text{ and therefore}$$

$$\det\left(\binom{j+k}{i-j+1}\right)_{i,j=0}^{n-1} = C_n^{(k)}. \tag{33}$$

**Remarks**

1) As a consequence of (22) we see (cf. [6]) that the inverse of $\left((-1)^{i-j}\binom{i+j}{i-j}\right)_{i,j \geq 0}$ is the "Catalan triangle"

$$\left((-1)^{i-j}\binom{i+j}{i-j}\right)^{-1} = \left(\frac{2j+1}{i+j+1}\binom{2i}{i-j}\right)_{i,j \geq 0} = \left(C_{i-j}^{(2j+1)}\right)_{i,j \geq 0}. \tag{34}$$

The first terms are

$$\begin{pmatrix} 1 & 0 & 0 & 0 & 0 & 0 & 0 \\ 1 & 1 & 0 & 0 & 0 & 0 & 0 \\ 2 & 3 & 1 & 0 & 0 & 0 & 0 \\ 5 & 9 & 5 & 1 & 0 & 0 & 0 \\ 14 & 28 & 20 & 7 & 1 & 0 & 0 \\ 42 & 90 & 75 & 35 & 9 & 1 & 0 \\ 132 & 297 & 275 & 154 & 54 & 11 & 1 \end{pmatrix}.$$

2) Let $A_n(x) = \left( \binom{x+i+j}{i-j+1} \right)_{i,j=0}^{n-1}$. The identity $\det A_n(x) = G_n(x) = \dfrac{x}{n!} \prod_{j=n+1}^{2n-1}(x+j)$

implies that

$$\det A_n(-m) = 0 \text{ for } n+1 \leq m \leq 2n-1, \tag{35}$$

if we set $\binom{x}{k} = 0$ for $k < 0$.

It is perhaps interesting to note that the corresponding null-spaces are closely related to the

Lucas polynomials $L_m(x) = \sum_{j=0}^{\lfloor m/2 \rfloor} (-1)^j \dfrac{m}{m-j} \binom{m-j}{j} x^{m-2j}$. More precisely let

$$v_{n,m} = \left( \dfrac{m}{m-j} \binom{m-j}{j} \right)_{0 \leq j \leq n-1}^t. \text{ Then}$$

$$A_n(-m) v_{n,m} = 0 \tag{36}$$

or equivalently

$$\sum_{j=0}^{n-1} \binom{i+j-m}{i+1-j} \dfrac{m}{m-j} \binom{m-j}{j} = 0 \tag{37}$$

for $0 \leq i \leq n-1$ and $n+1 \leq m \leq 2n-1$.

Identity (37) can be obtained from the Rothe-type formula (cf. [13],(5.62))

$$\sum_{j=0}^{n} \dfrac{r}{r+tj} \binom{r+tj}{j} \binom{s+tn-tj}{n-j} = \binom{r+s+tn}{n}. \tag{38}$$

If we set $r \to x$, $n \to i+1$, $t \to -1$, $s \to 2i+1-x$, we get

$$\sum_{j=0}^{i+1} \frac{x}{x-j} \binom{x-j}{j} \binom{i+j-x}{i-j+1} = \binom{i}{i+1} = 0 \text{ for all } i.$$

Thus (37) is true for $i < n-1$.

For $i = n-1$ we get

$$\sum_{j=0}^{n-1} \binom{i+j-x}{i+1-j} \frac{x}{x-j} \binom{x-j}{j} = -\binom{2n-1-x}{0} \frac{x}{x-n} \binom{x-n}{n} = \frac{x}{n!} \prod_{j=n+1}^{2n-1}(x-j)$$

because $m \leq 2n-1$.

The same argument also gives

$$\left( \binom{j-m}{i-j+1} \right)_{i,j=0}^{n-1} v_{n,m} = 0 \tag{39}$$

for $n+1 \leq m \leq 2n-1$.

Let me also mention the analogous results for Lucas polynomials and central binomial coefficients which are a sort of companion for the Fibonacci polynomials and Catalan numbers.

We consider a variant of the Lucas polynomials defined by

$$L_n(x) = \sum_{j=0}^{n} (-1)^j \frac{n}{n-j} \binom{n-j}{j} x^{n-2j} \tag{40}$$

for $n > 0$. and $L_0(x) = 1$. These polynomials are monic and orthogonal with $s(n) = 0$, $t(0) = 2$ and $t(n) = 1$ for $n > 0$. Their moments are $\Lambda(x^{2n}) = \binom{2n}{n}$ and $\Lambda(x^{2n+1}) = 0$.

More generally the numbers $c(n,k)$ are given by

$$c(2n+k, k) = \binom{2n+k}{n} \tag{41}$$

and $c(n,j) = 0$ else, because these numbers satisfy the identities (21) which reduce to

$$\binom{2n}{n} = 2\binom{2n-1}{n-1} \text{ and } \binom{2n+k}{n} = \binom{2n+k-1}{n} + \binom{2n+k-1}{n-1} \text{ for } k > 0.$$

Note that (41) is equivalent with

$$\sum_{k=0}^{\lfloor\frac{n}{2}\rfloor}\binom{n}{k}L_{n-2k}=x^{n}. \tag{42}$$

From $\Lambda\left(L_{2n}(x)\right)=\sum_{j=0}^{n}(-1)^{n-j}\frac{2n}{n+j}\binom{n+j}{n-j}\binom{2j}{j}=[n=0]$ we get

$$\det\left(\frac{2i+2}{i+j+1}\binom{i+j+1}{i-j+1}\right)_{i,j=0}^{n-1}=\binom{2n}{n}. \tag{43}$$

$$L_{2n+k}(x)=\sum_{0\leq 2j+k\leq 2n+k}(-1)^{n-j}\frac{2n+k}{n+k+j}\binom{n+k+j}{n-j}x^{2j+k}$$

implies $\dfrac{\Lambda\left(L_{2n+k}(x)L_k(x)\right)}{\Lambda\left(L_k(x)^2\right)}=\sum_{0\leq 2j+k\leq 2n+k}(-1)^{n-j}\dfrac{2n+k}{n+k+j}\binom{n+k+j}{n-j}\dfrac{\Lambda\left(x^{2j+k}L_k(x)\right)}{\Lambda\left(L_k(x)^2\right)}.$

Therefore we get by (23) and $\Lambda\left(x^{2j+k}L_k(x)\right)=0$ for $j<0$

$$\sum_{j=0}^{n}(-1)^{n-j}\frac{2n+k}{n+k+j}\binom{n+k+j}{n-j}\binom{2j+k}{j}=[n=0] \tag{44}$$

which implies

$$\det\left(\frac{2i+k+1}{i+j+k}\binom{i+j+k}{i-j+1}\right)_{i,j=0}^{n-1}=\binom{2n+k-1}{n}. \tag{45}$$

From

$$x^{n-1}L_{n+k+1}(x)=\sum_{0\leq 2j\leq n+k+1}(-1)^{j}\frac{n+k+1}{n+k+1-j}\binom{n+k+1-j}{j}x^{2n+k-2j}$$
$$=\sum_{j}(-1)^{n-j}\frac{n+k+1}{k+1+j}\binom{k+1+j}{n-j}x^{k+2j}$$

we get

$$\sum_{j}(-1)^{n-j}\frac{n+k+1}{k+1+j}\binom{k+1+j}{n-j}\binom{2j+k}{j}=\frac{\Lambda\left(x^{n-1}L_{n+k+1}(x)L_k(x)\right)}{\Lambda\left(L_k(x)^2\right)}=[n=0] \text{ and therefore}$$

$$\det\left(\frac{i+k+1}{j+k}\binom{j+k}{i-j+1}\right)_{i,j=0}^{n-1}=\binom{2n+k-1}{n}. \tag{46}$$

**Remark**

Let $B_n(x) = \left( \dfrac{2i+1+x}{i+j+x} \binom{i+j+x}{i-j+1} \right)_{i,j=0}^{n-1}$. Then (45) gives

$$\det B_n(x) = \binom{2n-1+x}{n} = \dfrac{\prod_{j=n}^{2n-1}(x+j)}{n!}.$$

This implies $\det B_n(-m) = 0$ for $n \leq m \leq 2n-1$. The corresponding eigenvectors are

$$u_{n,m} = \left( \binom{m-j}{j} \right)^t_{0 \leq i,j \leq n-1}, \text{ i.e.}$$

$$B_n(-m) \left( \binom{m-j}{j} \right)^t_{0 \leq i,j \leq n-1} = 0. \tag{47}$$

Note that $u_{n,m}$ is related to the coefficients of the Fibonacci polynomials $F_m(x)$.

We have to show that for $0 \leq i \leq n-1$ and $n \leq m \leq 2n-1$ we have

$$\sum_{j=0}^{n-1} \dfrac{2i+1-m}{i+j-m} \binom{i+j-m}{i+1-j} \binom{m}{m-j} = 0. \tag{48}$$

To prove this we show first that for all $i$

$$\sum_{j=0}^{i+1} \dfrac{2i+1-x}{i+j-x} \binom{i+j-x}{i-j+1} \binom{x-j}{j} = \sum_{j=0}^{i+1} \dfrac{2i+1-x}{2i+1-j-x} \binom{2i+1-j-x}{j} \binom{x-i-1+j}{i+1-j} = 0.$$

This follows from (38) if we let $r \to 2i+1-x$, $t \to -1$, $n \to i+1$, $s \to x$.

Thus (48) is true for all $m$ if $i < n-1$.

For $i = n-1$ we get

$$\sum_{j=0}^{n-1} \dfrac{2i+1-x}{i+j-x} \binom{i+j-x}{i-j+1} \binom{x-j}{j} = -\dfrac{2n-1-x}{2n-1-x} \binom{2n-1-x}{0} \binom{x-n}{n} = -\binom{x-n}{n}.$$

This vanishes for $n \leq x \leq 2n-1$.

The same argument also gives

$$\left( \dfrac{i+1-m}{j-m} \binom{i+j-m}{i-j+1} \right)_{i,j=0}^{n-1} \left( \binom{m-j}{j} \right)^t_{0 \leq i,j \leq n-1} = 0. \tag{49}$$

### 2.1.2. Expansion with respect to Gould polynomials

Consider the special Gould polynomials

$$G_n(x) = G_n(x,2) = \frac{x}{2n+x}\binom{2n+x}{n} = \frac{x}{n!}\prod_{j=n+1}^{2n-1}(x+j) \qquad (50)$$

for $n > 0$ and $G_0(x) = 1$. Let $Ef(x) = f(x+1)$ be the shift operator and $\Delta = E - 1$ the difference operator on the polynomials.

These Gould polynomials satisfy

$$E^{-2}\Delta G_n(x) = G_n(x-1) - G_n(x-2) = G_{n-1}(x) \qquad (51)$$

because $G_n(x)$ is a polynomial of degree $n$ and (30) shows that (51) is true for infinitely many non-negative integers $x = k$.

Let us expand the polynomial $\binom{x-1}{n}$ as linear combinations of Gould polynomials,

$$\binom{x-1}{n} = \sum_{j=0}^{n} b(n,j) G_j(x).$$

Since $\left(E^{-2}\Delta\right)^m G_n(x)\Big|_{x=0} = G_{n-m}(x)\Big|_{x=0} = [n = m]$ we get

$$b(n,j) = \left(E^{-2}\Delta\right)^j \binom{x-1}{n}\bigg|_{x=0} = \binom{x-2j-1}{n-j}\bigg|_{x=0} = \binom{-2j-1}{n-j} = (-1)^{n-j}\binom{n+j}{n-j}$$

Thus we have

$$\binom{x-1}{n} = \sum_{j=0}^{n}(-1)^{n-j}\binom{n+j}{n-j}G_j(x). \qquad (52)$$

For $x = 1$ this gives (2) and thus also (1).

If we expand $\binom{x-k}{n} = \sum_{j=0}^{n} b(n,j) G_j(x)$ we get in the same way

$$\binom{x-k}{n} = \sum_{j=0}^{n}(-1)^{n-j}\binom{n+j+k-1}{n-j}G_j(x). \qquad (53)$$

For $x = k$ this gives $\binom{0}{n} = \sum_{j=0}^{n}(-1)^{n-j}\binom{n+j+k-1}{n-j}C_j^{(k)},$

which implies

$$\det\left(\binom{i+j+k}{i-j+1}\right)_{i,j=0}^{n-1} = C_n^{(k)}. \tag{54}$$

Expanding the polynomials $\binom{x+n-1}{n}$ we get

$\binom{x+n-1-k}{n} = \sum_{j=0}^{n}(-1)^{n-j}\binom{j+k}{n-j}G_j(x)$ which for $x = k$ reduces to

$\sum_{j=0}^{n}(-1)^{n-j}\binom{j+k}{n-j}C_j^{(k)} = \binom{n-1}{n} = [n=0].$ This implies

$$\det\left(\binom{j+k}{i-j+1}\right)_{i,j=0}^{n-1} = C_n^{(k)}. \tag{55}$$

**Remark**

Formulae (54) and (55) are equivalent by elementary row operations using Vandermonde's identity

$$\sum_{\ell=0}^{i}\binom{i}{i-\ell}\binom{j+k}{\ell-j} = \binom{i+j+k}{i-j}.$$

Formula (53) can also be obtained from (38). If we set $r = x$, $t = 2$, $s = -k - 2n$ we get (53) because

$$\sum_{j=0}^{n}\frac{x}{x+2j}\binom{x+2j}{j}\binom{-k-2j}{n-j} = \sum_{j=0}^{n}\frac{x}{x+2j}\binom{x+2j}{j}(-1)^{n-j}\binom{n+j+k-1}{n-j} = \binom{x-k}{n}.$$

The general Gould polynomials $G_n(x,r) = \frac{x}{rn+x}\binom{rn+x}{n}$ satisfy

$E^{-r}\Delta G_n(x,r) = G_{n-1}(x,r)$ (cf. [17], p.55 ). In the same way as above we get

$$\binom{x-k}{n} = \sum_{j=0}^{n} (-1)^{n-j} \binom{n+(r-1)j+k-1}{n-j} G_j(x,r) \tag{56}$$

which for $x = k$ reduces to

$$\sum_{j=0}^{n} (-1)^{n-j} \binom{n+(r-1)j+k-1}{n-j} \frac{k}{k+rj} \binom{k+rj}{j} = [n=0]. \tag{57}$$

This implies

$$\det\left(\binom{i+(r-1)j+k}{i-j+1}\right)_{i,j=0}^{n-1} = \frac{k}{rn+k} \binom{rn+k}{n}. \tag{58}$$

If we set $k = \alpha + 1 - n$ in (56) we get

$$\binom{x-k}{n} = \binom{x+n-\alpha-1}{n} = (-1)^n \binom{\alpha-x}{n} \quad \text{and therefore}$$

$$\sum_{j=0}^{n} (-1)^{n-j} \binom{(r-1)j+\alpha}{n-j} \frac{\gamma}{rj+\gamma} \binom{rj+\gamma}{j} = (-1)^n \binom{\alpha-\gamma}{n}. \tag{59}$$

A combinatorial proof of this identity has been given by Ricky X.F. Chen [3].

For $\alpha = \gamma = k$ we get

$$\sum_{j=0}^{n} (-1)^{n-j} \binom{(r-1)j+k}{n-j} \frac{k}{k+rj} \binom{k+rj}{j} = [n=0] \tag{60}$$

and therefore

$$\det\left(\binom{(r-1)j+k}{i-j+1}\right)_{i,j=0}^{n-1} = \frac{k}{rn+k} \binom{rn+k}{n}. \tag{61}$$

The last identity has also been obtained in [11].

### 2.1.3. Proof by means of hypergeometric identities

**2.1.3.1.** Using the notation $(x)_n = x(x+1)\cdots(x+n-1)$ and Chu-Vandermonde's formula (cf. [1], 2.2.3)

$$\sum_{j=0}^{n} \frac{(a)_j (-n)_j}{j!(c)_j} = {}_2F_1\left(\begin{array}{c} a, -n \\ c \end{array}; 1\right) = \frac{(a-c)_n}{(c)_n} \tag{62}$$

we get another proof of (2).

After changing $j \to n-j$ we get

$$\sum_{j=0}^{n}(-1)^{n-j}\binom{n+j}{n-j}C_j = \sum_{j=0}^{n}(-1)^j\binom{2n-j}{j}\frac{1}{(n-j+1)}\binom{2n-2j}{n-j}$$

$$= \sum_{j=0}^{n}(-1)^j \frac{(2n-j)!}{j!(2n-2j)!}\frac{(2n-2j)!}{(n-j)!(n+1-j)!} = \sum_{j=0}^{n}\frac{(2n)!}{j!(-2n)_j}\frac{(-n)_j}{n!}\frac{(-n-1)_j}{(n+1)!}$$

$$= C_n \sum_{j=0}^{n}\frac{(-n)_j(-n-1)_j}{j!(-2n)_j} = C_n \,{}_2F_1\!\left[\begin{matrix}-n-1,-n\\-2n\end{matrix};1\right] = C_n \frac{(1-n)_n}{(-2n)_n} = [n=0].$$

**2.1.3.2.** After completion of a first draft of this paper Darij Grinberg showed me a reference to an elementary short proof by Max Alexeyev (cf.
https://artofproblemsolving.com/community/c6h349513p1878114): Writing

$$C_k\binom{n+k}{n-k} = \frac{(2k)!}{k!(k+1)!}\frac{(n+k)!}{(n-k)!(2k)!} = \frac{(n+k)!}{k!n!}\frac{1}{(n+1)}\frac{(n+1)!}{(k+1)!(n-k)!} = \frac{1}{(n+1)}\binom{n+k}{k}\binom{n+1}{n-k}$$

and using $(-1)^k\binom{n+k}{k} = \binom{-n-1}{k}$ we get

$$\sum_{j=0}^{n}(-1)^{n-j}\binom{n+j}{n-j}C_j = \frac{(-1)^n}{n+1}\sum_{k=0}^{n}\binom{-n-1}{k}\binom{n+1}{n-k} = \frac{(-1)^n}{n+1}\binom{0}{n} = [n=0]$$

by Vandermonde's identity.

The same method also gives (31):

$$\sum_{j=0}^{n}(-1)^{n-j}\binom{n+k+j}{n-j}C_j^{(k+1)} = (-1)^n\sum_{j=0}^{n}(-1)^j \frac{(n+k+j)!}{(n-j)!(2j+k)!}\frac{(k+1)}{(2j+k+1)}\binom{2j+k+1}{j}$$

$$= (-1)^n\sum_{j=0}^{n}(-1)^j \frac{(n+k+j)!}{(n-j)!}\frac{(k+1)}{j!(j+k+1)!} = (-1)^n \frac{(k+1)}{n+k+1}\sum_{j=0}^{n}(-1)^j\binom{n+k+j}{j}\binom{n+k+1}{n-j}$$

$$= (-1)^n \frac{(k+1)}{n+k+1}\sum_{j=0}^{n}\binom{-n-k-1}{j}\binom{n+k+1}{n-j} = [n=0].$$

**2.2.** A different method to prove (1) uses

**Lemma 3 (C. Krattenthaler [15], Th. 27)**

*With the usual notations of $q$-calculus we have*

$$\det\!\left(q^{jL_i}\begin{bmatrix}L_i+A-j\\L_i+j\end{bmatrix}\right) = q^{\sum_{i=1}^{n}iL_i}\prod_{i=1}^{n}\frac{[L_i+A-n]!}{[L_i+n]![A-2i]!}\prod_{j=1}^{n}\prod_{i=1}^{j-1}[L_i-L_j][L_i+L_j+A+1]. \quad (63)$$

For $q = 1$ this reduces to

$$\det\left(\binom{L_i + A - j}{L_i + j}\right)_{i,j=1}^n = \prod_{i=1}^n \frac{(L_i + A - n)!}{(L_i + n)!(A - 2i)!} \prod_{j=1}^n \prod_{i=1}^{j-1} (L_i - L_j)(L_i + L_j + A + 1). \tag{64}$$

In order to apply this to $\det\left(\binom{i+j+1}{i-j+1}\right)_{i,j=0}^{n-1}$ we first change $i \to i-1$ and $j \to j-1$ and then reverse the order of the rows and columns in the matrix $\left(\binom{i+j-1}{i-j+1}\right)_{i,j=1}^n$. This means we change $i \to n+1-i$ and $j \to n+1-j$.

We then get the matrix $\left(\binom{2n+1-i-j}{j-i+1}\right)_{i,j=1}^n$.

Since this operation lets the determinant unchanged we get

$\det A_n = \det\left(\binom{2n+1-i-j}{j-i+1}\right)_{i,j=1}^n$. Choosing $L_i = 1-i$ and $A = 2n$ in (64) we get

$$\det A_n = \prod_{i=1}^n \frac{(1-i+n)!}{(1-i+n)!(2n-2i)!} \prod_{j=1}^n \prod_{i=1}^{j-1} (j-i)(2-i-j+2n+1)$$

$$= \prod_{i=1}^n \frac{(i-1)!}{(2n-2i)!} \prod_{j=1}^n \prod_{i=1}^{j-1} (3-i-j+2n).$$

In order to simplify this let $f(n) = \prod_{i=1}^n \frac{(i-1)!}{(2n-2i)!} = \prod_{i=1}^n \frac{(i-1)!}{(2i-2)!} = \frac{(n-1)!}{(2n-2)!} f(n-1)$

and

$$g(n) = \prod_{j=1}^n \prod_{i=1}^{j-1} (3-i-j+2n) = \prod_{j=1}^n \frac{(2n+2-j)!}{(2n+3-2j)!}$$

$$= \frac{(2n)!}{(2n-1)!} \frac{(2n-1)!}{(2n-3)!} \cdots \frac{(n+2)!}{(3)!} = \frac{(2n)!}{(n+1)!} g(n-1).$$

This gives

$$f(n)g(n) = \frac{(n-1)!}{(2n-2)!} \frac{(2n)!}{(n+1)!} f(n-1)g(n-1) = \frac{(2n)(2n-1)}{n(n+1)} f(n-1)g(n-1) = \frac{1}{n+1}\binom{2n}{n}.$$

The same method also leads to

**Theorem 4**

$$\det\left(\binom{i+j+m}{i-j+m}\right)_{i,j=0}^{n-1} = \prod_{j=1}^{n-1}\frac{j!}{(2j)!}\prod_{j=1}^{n-1}\frac{(2m+2j)!}{(2m+j)!} = \det\left(C_{n+i+j}\right)_{i,j=0}^{m-1}. \qquad (65)$$

**Proof**

First we write $\det\left(\binom{i+j+m}{i-j+m}\right)_{i,j=0}^{n-1} = \det\left(\binom{i+j+m-2}{i-j+m}\right)_{i,j=1}^{n}$ and then we again reverse the order of the rows and columns. Comparing with (64) we choose $L_i = m-i$ and $A = 2n$ and get

$$\det\left(\binom{i+j+m}{i-j+m}\right)_{i,j=0}^{n-1} = \det\left(\binom{2n+m-i-j}{j-i+m}\right)_{i,j=1}^{n}$$

$$= \prod_{i=1}^{n}\frac{(m-i+n)!}{(m-i+n)!(2n-2i)!}\prod_{j=1}^{n}\prod_{i=1}^{j-1}(j-i)(2m+1-i-j+2n)$$

$$= \prod_{i=1}^{n}\frac{(i-1)!}{(2n-2i)!}\prod_{j=1}^{n}\frac{(2m-j+2n)!}{(2m+1+2n-2j)!} = f(n)g(n,m)$$

with $f(n) = \prod_{i=1}^{n}\frac{(i-1)!}{(2n-2i)!}$ and

$$g(n,m) = \prod_{j=1}^{n}\frac{(2m-j+2n)!}{(2m+1+2n-2j)!} = \frac{(2m+2n-1)!(2m+2n-2)!\cdots(2m+n)!}{(2m+2n-1)!(2m+2n-3)!\cdots(2m+1)!}.$$

This gives $\dfrac{f(n)}{f(n-1)} = \dfrac{(n-1)!}{(2n-2)!}$ and $\dfrac{g(n,m)}{g(n-1,m)} = \dfrac{(2m+2n-2)!}{(2m+n-1)!}$

because

$$g(n-1,m) = \frac{(2m+2n-3)!(2m+2n-4)!\cdots(2m+n-1)!}{(2m+2n-3)!(2m+2n-5)!\cdots(2m+1)!}.$$

Therefore we have

$$\frac{f(n)g(n,m)}{f(n-1)g(n-1,m)} = \frac{(n-1)!}{(2n-2)!}\frac{(2m+2n-2)!}{(2m+n-1)!} \qquad (66)$$

On the other hand it is well known that (cf.e.g. [16],Theorem 33)

$$\det\left(C_{n+i+j}\right)_{i,j=0}^{m-1} = \prod_{j=1}^{n-1}\prod_{i=1}^{j}\frac{2m+i+j}{i+j} = \prod_{j=1}^{n-1}\frac{j!}{(2j)!}\prod_{j=1}^{n-1}\frac{(2m+2j)!}{(2m+j)!} \qquad (67)$$

Comparing (66) with (67) we get (65). $\square$

**Remark**

As already mentioned the appearance of the Hankel determinants $\det\left(C_{n+i+j}\right)_{i,j=0}^{m-1}$ led to the conjecture that this is a special case of a general theorem about Hankel determinants which in the mean-time has been proved by M. Tyson [18]:

**Theorem 5**

Let $p_n(x) = \sum_{j=0}^{n} (-1)^{n-j} p(n,j) x^j$ be monic orthogonal polynomials with moments $M_n$. Then

$$\det\left(p(i+m,j)\right)_{i,j=0}^{n-1} = \frac{\det\left(M_{n+i+j}\right)_{i,j=0}^{m-1}}{\det\left(M_{i+j}\right)_{i,j=0}^{m-1}}. \tag{68}$$

We shall not reproduce the proof, but state only some consequences.

Formula (68) generalizes the formulas

$$\det\left(M_{i+j+1}\right)_{i,j=0}^{m-1} = (-1)^m p_m(0) \det\left(M_{i+j}\right)_{i,j=0}^{m-1} \tag{69}$$

and

$$\det\left(M_{i+j+2}\right)_{i,j=0}^{m-1} = v(m) \det\left(M_{i+j}\right)_{i,j=0}^{m-1} \tag{70}$$

with $v(m) = \sum_{k=0}^{m} p_k(0)^2 t(k) t(k+1) \cdots t(m-1)$ (cf. [16],(5.41)) if
$p_n(x) = (x - s(n-1)) p_{n-1}(x) - t(n-2) p_{n-2}(x)$.

Let us deduce (70) from (68): We have

$p_n(0) = -s(n-1) p_{n-1}(0) - t(n-2) p_{n-2}(0)$ and
$p'_n(0) = -s(n-1) p'_{n-1}(0) - t(n-2) p'_{n-2}(0) + p_{n-1}(0)$.

Therefore we get

$$v(m) = \det\begin{pmatrix} p(m,0) & p(m,1) \\ p(m+1,0) & p(m+1,1) \end{pmatrix}$$

$$= -\det\begin{pmatrix} p_m(0) & p'_m(0) \\ s(m) p_m(0) - t(m-1) p_{m-1}(0) & s(m) p'_m(0) + t(m-1) p'_{m-1}(0) - p_m(0) \end{pmatrix}$$

$$= -\det\begin{pmatrix} p_m(0) & p'_m(0) \\ -t(m-1) p_{m-1}(0) & t(m-1) p'_{m-1}(0) - p_m(0) \end{pmatrix}$$

$$= p_m(0)^2 + t(m-1)p_{m-1}(0)p'_m(0) - t(m-1)p'_{m-1}(0)p_m(0)$$
$$= p_m(0)^2 + t(m-1)\det\begin{pmatrix} p(m-1,0) & p(m-1,1) \\ p(m,0) & p(m,1) \end{pmatrix} = p_m(0)^2 + t(m-1)v(m-1).$$

Since $v(0) = 1$ we get $v(m) = \sum_{k=0}^{m} p_k(0)^2 t(k)t(k+1)\cdots t(m-1)$.

As another simple application consider again the polynomials
$$p_n(x) = \sum_{j=0}^{n} (-1)^{n-j} \begin{bmatrix} n \\ j \end{bmatrix} q^{(n-1)(n-j)} x^j \text{ with moments } M_n = q^{\binom{n}{2}}. \text{ Here we get}$$

$$\det\left( q^{(i+m-1)(i+m-j)} \begin{bmatrix} i+m \\ j \end{bmatrix} \right)_{i,j=0}^{n-1} = \frac{\det\left( q^{\binom{n+i+j}{2}} \right)_{i,j=0}^{m-1}}{\det\left( q^{\binom{i+j}{2}} \right)_{i,j=0}^{m-1}},$$

which can be simplified to

$$\det(A_n(m)) = \det\left( q^{\binom{i-j}{2}} \begin{bmatrix} i+m \\ j \end{bmatrix} \right)_{i,j=0}^{n-1} = 1 \tag{71}$$

since

$$\det\left( q^{(i+m-1)(i+m-j)} \begin{bmatrix} i+m \\ j \end{bmatrix} \right)_{i,j=0}^{n-1} = \det\left( q^{\binom{i}{2} - \binom{j}{2} + m^2 - m + 2im - jm} q^{\binom{i-j}{2}} \begin{bmatrix} i+m \\ j \end{bmatrix} \right)_{i,j=0}^{n-1}$$
$$= q^{2\binom{m}{2}n + m\binom{n}{2}} \det\left( q^{\binom{i-j}{2}} \begin{bmatrix} i+m \\ j \end{bmatrix} \right)_{i,j=0}^{n-1}$$

and

$$\frac{\det\left( q^{\binom{n+i+j}{2}} \right)_{i,j=0}^{m-1}}{\det\left( q^{\binom{i+j}{2}} \right)_{i,j=0}^{m-1}} = \frac{\det\left( q^{\binom{n}{2} + in + jn + \binom{i+j}{2}} \right)_{i,j=0}^{m-1}}{\det\left( q^{\binom{i+j}{2}} \right)_{i,j=0}^{m-1}} = q^{\binom{n}{2}m + 2n\binom{m}{2}}.$$

**Remark**

Of course there are simpler proofs of (71). Let for example $B_n(m) = \left(b_m(i,j)\right)_{i,j=0}^{n-1}$ be the upper triangular matrix satisfying $b_m(i,i) = 1$ and $b_m(i,i+1) = q^m$ for all $i$ and $b_m(i,j) = 0$ else. Then $\det B_n(m) = 1$ and $A_n(m) = A_n(m-1)B_n(m)$ because

$$q^m q^{\binom{i-(j-1)}{2}} \begin{bmatrix} m-1+i \\ j-1 \end{bmatrix} + q^{\binom{i-j}{2}} \begin{bmatrix} m-1+i \\ j \end{bmatrix} = q^{\binom{i-j}{2}} \begin{bmatrix} m+i \\ j \end{bmatrix}. \text{ Thus}$$

$$\det A_n(m) = \det A_n(m-1) = \cdots = \det A_n(0) = \det\left(q^{\binom{i-j}{2}} \begin{bmatrix} i \\ j \end{bmatrix}\right)_{i,j=0}^{n-1} = 1.$$

For another illustration of Theorem 5 consider the sequence $a(n) = \dfrac{1}{n+1}$, which gives the famous Hilbert matrix.

The corresponding orthogonal polynomials are (cf. [9])

$$r_n(x) = \sum_{j=0}^{n} (-1)^{n-j} \frac{\binom{n}{j}\binom{n+j}{j}}{\binom{2n}{n}} x^j.$$

Thus (68) reduces to

$$\det\left(\frac{\binom{i+m}{j}\binom{i+m+j}{j}}{\binom{2i+2m}{i+m}}\right)_{i,j=0}^{n-1} = \frac{\det\left(\dfrac{1}{n+i+j+1}\right)_{i,j=0}^{m-1}}{\det\left(\dfrac{1}{i+j+1}\right)_{i,j=0}^{m-1}}. \tag{72}$$

Let us verify this identity by direct computation:

To compute the determinant of $\left(\dfrac{1}{i+j+n+1}\right)_{i,j=0}^{m-1} = \left(\dfrac{1}{i+j+n-1}\right)_{i,j=1}^{m}$ we use Cauchy's formula (cf. [15],(2.7))

$$\det\left(\frac{1}{x_i+y_j}\right)_{i,j=1}^{m} = \frac{\prod_{1 \leq i < j \leq m}(x_i-x_j)(y_i-y_j)}{\prod_{1 \leq i,j \leq m}(x_i+y_j)}$$

with $x_i = i+n-1, y_j = j$. This gives

$$\det\left(\frac{1}{n+i+j+1}\right)_{i,j=0}^{m-1} = \prod_{j=0}^{m-1} \frac{j!\,j!(n+j)!}{(n+m+j)!} \qquad (73)$$

and thus

$$\frac{\det\left(\frac{1}{n+i+j+1}\right)_{i,j=0}^{m-1}}{\det\left(\frac{1}{i+j+1}\right)_{i,j=0}^{m-1}} = \frac{\prod_{j=0}^{m-1}\frac{j!\,j!(n+j)!}{(n+m+j)!}}{\prod_{j=0}^{m-1}\frac{j!\,j!\,j!}{(m+j)!}} = \prod_{j=0}^{m-1}\frac{(n+j)!(m+j)!}{(n+m+j)!\,j!}.$$

For the left-hand side of (72) we get using (65)

$$\det\left(\frac{\binom{i+m}{j}\binom{i+m+j}{j}}{\binom{2i+2m}{i+m}}\right)_{i,j=0}^{n-1} = \det\left(\frac{(i+m)!^2(2j)!}{j!\,j!(2i+2m)!}\binom{i+j+m}{i-j+m}\right)_{i,j=0}^{n-1}$$

$$= \det\left(\binom{i+j+m}{i-j+m}\right)_{i,j=0}^{n-1} \prod_{j=0}^{n-1}\frac{(i+m)!^2(2j)!}{j!\,j!(2i+2m)!}$$

$$= \prod_{j=1}^{n-1}\frac{j!}{(2j)!}\frac{(2m+2j)!}{(2m+j)!}\frac{(2j)!}{j!\,j!}\frac{(j+m)!(j+m)!}{(2j+2m)!} = \prod_{j=1}^{n-1}\frac{(j+m)!(j+m)!}{j!(2m+j)!}.$$

To verify (72) we must show that

$$u(n,m) = \prod_{j=1}^{n-1}\frac{(j+m)!(j+m)!}{j!(2m+j)!} = v(n,m) = \prod_{j=0}^{m-1}\frac{(n+j)!(m+j)!}{j!(n+m+j)!}.$$

This follows from $u(0,m) = v(0,m)$ and $\dfrac{u(n,m)}{u(n-1,m)} = \dfrac{(n+m-1)!(n+m-1)!}{(n-1)!(2m+n-1)!}$ and

$$\frac{v(n,m)}{v(n-1,m)} = \prod_{j=0}^{m-1}\frac{(n+j)}{(n+m+j)} = \frac{(n+m-1)!}{(n-1)!}\frac{(n+m-1)!}{(n+2m-1)!}.$$

As a further generalization we get

**Theorem 6**

$$\det\left(\binom{i+j+k+m}{i-j+m}\right)_{i,j=0}^{n-1} = \det\left(\binom{2n+m+k-i-j}{j-i+m}\right)_{i,j=1}^{n} = \det\left(C^{(2i+k+1)}_{n-i+j}\right)_{i,j=0}^{m-1}. \qquad (74)$$

Let us first show that this indeed generalizes Theorem 4.

First we see that

$$\det\begin{pmatrix} C_n & C_{n+1} \\ C_{n+1} & C_{n+2} \end{pmatrix} = \det\begin{pmatrix} C_n^{(1)} & C_{n+1}^{(1)} \\ C_{n-1}^{(3)} & C_n^{(3)} \end{pmatrix}$$

From $c(z)^2 = c(z) + zc(z)^3$ and $C_n^{(2)} = C_{n+1}^{(1)} = C_{n+1}$ we get $C_{n+1} = C_n + C_{n-1}^{(3)}$.

Thus $\det\begin{pmatrix} C_n & C_{n+1} \\ C_{n+1} & C_{n+2} \end{pmatrix} = \det\begin{pmatrix} C_n & C_{n+1} \\ C_{n+1} - C_n & C_{n+2} - C_{n+1} \end{pmatrix} = \det\begin{pmatrix} C_n^{(1)} & C_{n+1}^{(1)} \\ C_{n-1}^{(3)} & C_n^{(3)} \end{pmatrix}.$

In the same way we get

$$\det\begin{pmatrix} C_n & C_{n+1} & C_{n+2} \\ C_{n+1} & C_{n+2} & C_{n+3} \\ C_{n+2} & C_{n+3} & C_{n+4} \end{pmatrix} = \det\begin{pmatrix} C_n & C_{n+1} & C_{n+2} \\ C_{n+1} & C_{n+2} & C_{n+3} \\ C_n^{(3)} & C_{n+1}^{(3)} & C_{n+2}^{(3)} \end{pmatrix} = \det\begin{pmatrix} C_n & C_{n+1} & C_{n+2} \\ C_{n-1}^{(3)} & C_n^{(3)} & C_{n+1}^{(3)} \\ C_n^{(3)} & C_{n+1}^{(3)} & C_{n+2}^{(3)} \end{pmatrix}$$

$$= \det\begin{pmatrix} C_n & C_{n+1} & C_{n+2} \\ C_{n-1}^{(3)} & C_n^{(3)} & C_{n+1}^{(3)} \\ C_{n-2}^{(5)} & C_{n-1}^{(5)} & C_n^{(5)} \end{pmatrix} = \det\begin{pmatrix} C_n^{(1)} & C_{n+1}^{(1)} & C_{n+2}^{(1)} \\ C_{n-1}^{(3)} & C_n^{(3)} & C_{n+1}^{(3)} \\ C_{n-2}^{(5)} & C_{n-1}^{(5)} & C_n^{(5)} \end{pmatrix}.$$

Analogously in the general case.

**Proof of Theorem 6**

Comparing with (64) we choose $L_i = -i + m$ and $A = 2n + k$. Then we get

$$w(n,m,k) = \det\left(\binom{2n+m+k-i-j}{j-i+m}\right)_{i,j=1}^n$$

$$= \prod_{i=1}^n \frac{(k+m-i+n)!}{(m-i+n)!(2n+k-2i)!} \prod_{j=1}^n \prod_{i=1}^{j-1}(j-i)(2m+1-i-j+2n+k)$$

$$= \prod_{i=1}^n \frac{(i-1)!(k+m-i+n)!}{(2n+k-2i)!(m-i+n)!} \prod_{j=1}^n \frac{(2m-j+2n+k)!}{(2m+1+2n-2j+k)!}.$$

This implies

$$v(n,m,k) = \frac{w(n,m,k)}{w(n-1,m,k)} = \frac{(n-1)!(k+m-1+n)!(2m+2n+k-2)!}{(2n+k-2)!(m+n-1)!(n+2m+k-1)!}$$

$$= \frac{\prod_{\ell=0}^{2m-1}(2n-1+k+\ell)}{\prod_{\ell=0}^{m-1}(n+\ell)(n+k+\ell+m)}. \tag{75}$$

for

$$\frac{w(n,m,k)}{w(n-1,m,k)} = \frac{\prod_{i=1}^{n} \frac{(i-1)!(k+m-i+n)!}{(2n+k-2i)!(m-i+n)!} \prod_{j=1}^{n} \frac{(2m-j+2n+k)!}{(2m+1+2n-2j+k)!}}{\prod_{i=1}^{n-1} \frac{(i-1)!(k+m-i+n-1)!}{(2n+k-2-2i)!(m-i+n-1)!} \prod_{j=1}^{n-1} \frac{(2m-j+2n+k-2)!}{(2m-1+2n-2j+k)!}}$$

$$= \frac{(n-1)!(k+m-1+n)!(2m+2n+k-2)!}{(2n+k-2)!(m+n-1)!(n+2m+k-1)!} = \frac{1}{\prod_{\ell=0}^{m-1}(n+\ell) \prod_{\ell=0}^{m-1}(n+m+\ell+k)} \prod_{\ell=0}^{2m-1}(2n-1+k+\ell)$$

Thus

$$\det\left(\binom{2n+m+k-i-j}{j-i+m}\right)_{i,j=1}^{n} = \prod_{j=1}^{n} v(j,m,k).$$

Consider now the right-hand side $M(m,n,k) = \det\left(C_{n-i+j}^{(k+1+2i)}\right)_{i,j=0}^{m-1}$.

This determinant can be computed with the Condensation method (cf. [15], Prop.10). This gives

$$M(m,n,k)M(m-2,n,k+2) \\ = M(m-1,n,k+2)M(m-1,n,k) - M(m-1,n+1,k)M(m-1,n-1,k+2). \quad (76)$$

Since no determinant vanishes this determines $M(m,n,k)$ by induction with $m$.

Let us verify the first terms. We have $M(0,n,k) = 1 = \prod_{j=1}^{n} v(j,0,k)$,

$$\prod_{j=1}^{n} v(j,1,k) = \prod_{j=1}^{n} \frac{(2j-1+k)(2j+k)}{j(j+1+k)} = \frac{1+k}{1+n+k}\binom{2n+k}{n} = C_n^{(k+1)}$$

Since

$$\frac{v(j,m,k)v(j,m-2,k+2)}{v(j,m-1,k)v(j,m-1,k+2)} = \frac{(k+j+m)(j+m-1)}{(k+j+m-1)(j+m-2)}$$

we get

$$\prod_{j=1}^{n-1} \frac{v(j,m,k)v(j,m-2,k+2)}{v(j,m-1,k)v(j,m-1,k+2)} = \frac{(m-1)(k+m)}{(n+m-2)(m+k+n-1)}$$

Therefore (76) reduces to the easily verified identity

$$\frac{(m-1)(m+k)}{(n+m-2)(m+k+n-1)} v(n,m,k)v(n,m-2,k+2) \\ = v(n,m-1,k+2)v(n,m-1,k) - v(n+1,m-1,k)v(n,m-1,k).$$

## 3. q-analogues

### 3.1. The Carlitz q-Catalan numbers

The Carlitz $q$ – Catalan numbers $c_n(q)$ satisfy $c_n(q) = \sum_{k=0}^{n-1} q^k c_k(q) c_{n-1-k}(q)$ with $c_0(q) = 1$.

The corresponding orthogonal polynomials are given by $f_n(x,q) = \sum_{k=0}^{\left\lfloor \frac{n}{2} \right\rfloor} (-1)^k q^{2\binom{k}{2}} \begin{bmatrix} n-k \\ k \end{bmatrix} x^{n-2k}$.

Here we have $s(n) = 0$ and $t(n) = q^n$ (cf. e.g.[5]). In this case the method of orthogonal polynomials gives nice $q$ – analogues, but unfortunately there are no closed formulas for $c_n(q)$ and the corresponding $c_n^{(k)}(q)$.

Let me only mention the analogue of (1).

$$\det\left( q^{2\binom{i-j}{2}} \begin{bmatrix} i+1+j \\ i+1-j \end{bmatrix} \right)_{i,j=0}^{n-1} = c_n(q). \tag{77}$$

The first terms are

$$\det(1) = 1, \quad \det\begin{pmatrix} 1 & q^2 \\ 1 & [3] \end{pmatrix} = 1 + q = c_2(q),$$

$$\det\begin{pmatrix} 1 & q^2 & 0 \\ 1 & [3] & q^2 \\ q^2 & (1+q^2)[3] & [5] \end{pmatrix} = 1 + 2q + q^2 + q^3 = c_3(q), \cdots.$$

**Remark**

From

$$c_n\left(\frac{1}{q}\right) = \det\left( \left(\frac{1}{q}\right)^{2\binom{i-j}{2}} \begin{bmatrix} i+j+1 \\ i-j+1 \end{bmatrix}_{\frac{1}{q}} \right)_{i,j=0}^{n-1} = \det\left( q^{-2\binom{i-j}{2} - 2j - 2ij + 2j^2} \begin{bmatrix} i+j+1 \\ i-j+1 \end{bmatrix}_q \right)_{i,j=0}^{n-1}$$

$$= \det\left( q^{-i^2 + i - 3j + j^2} \begin{bmatrix} i+j+1 \\ i-j+1 \end{bmatrix}_q \right)_{i,j=0}^{n-1} = q^{-2\binom{n}{2}} \det\left( \begin{bmatrix} i+j+1 \\ i-j+1 \end{bmatrix}_q \right)_{i,j=0}^{n-1}$$

we also get

$$\det\left( \begin{bmatrix} i+1+j \\ i+1-j \end{bmatrix} \right)_{i,j=0}^{n-1} = q^{2\binom{n}{2}} c_n\left(\frac{1}{q}\right). \tag{78}$$

It is perhaps interesting to consider the case $q=-1$. Let us first note that $\begin{bmatrix}2n\\2k\end{bmatrix}_{-1}=\binom{n}{k}$,

$$\begin{bmatrix}2n\\2k+1\end{bmatrix}_{-1}=0,\ \begin{bmatrix}2n+1\\2k\end{bmatrix}_{-1}=\binom{n}{k},\ \begin{bmatrix}2n+1\\2k+1\end{bmatrix}_{-1}=\binom{n}{k}.$$

This follows from $\begin{bmatrix}n\\k\end{bmatrix}_q = q^k\begin{bmatrix}n-1\\k\end{bmatrix}_q + \begin{bmatrix}n-1\\k-1\end{bmatrix}_q$ by induction with respect to $n$.

Further we get $c_{2n}(-1)=[n=0]$ and $c_{2n+1}(-1)=(-1)^n C_n$. This also follows by induction using $c_n(-1)=\sum_{k=0}^{n-1}(-1)^k c_k(-1) c_{n-1-k}(-1)$. For

$$c_{2n}(-1)=\sum_{k=0}^{2n-1}(-1)^k c_k(-1)c_{2n-1-k}(-1)=c_0(-1)c_{2n-1}(-1)+(-1)^{2n-1}c_{2n-1}(-1)c_0(-1)=0$$

and $c_{2n+1}(-1)=\sum_{k=0}^{2n}(-1)^k c_k(-1)c_{2n-k}(-1)=-\sum_{0\le 2j+1\le 2n-1}(-1)^j C_j(-1)^{n-1-j}C_{n-1-j}=(-1)^n C_n$.

Therefore (77) gives another class of matrices related to Catalan numbers:

$$\det\left(\left(\begin{bmatrix}i+j+1\\2\\j\end{bmatrix}\right)\right)_{i,j=0}^{2n}=(-1)^n C_n. \qquad (79)$$

In [4] we have introduced a $q-$analogue $g_n(r,q)$ of $\dfrac{1}{rn+1}\binom{rn+1}{n}$ defined by

$$g_n(r,q)=\sum_{k_1+\cdots+k_r=n-1}\prod_{j=1}^r\left(q^{(r-j)k_j}g_{k_j}(r,q)\right)\text{ with initial value }g_0(r,q)=1.$$

[4],(10), gives a combinatorial proof of the identity

$$\sum_{j=0}^n (-1)^{n-j} q^{r\binom{n-j}{2}}\begin{bmatrix}(r-1)j+1\\n-j\end{bmatrix}g_j(r,q)=[n=0]. \qquad (80)$$

This implies

$$\det\left(q^{r\binom{i-j+1}{2}}\begin{bmatrix}(r-1)j+1\\i-j+1\end{bmatrix}\right)_{i,j=0}^{n-1}=g_n(r,q), \qquad (81)$$

which generalizes (77) and (61) for $k=1$.

## 3.2. The q-Catalan numbers $C_n(q) = \frac{1}{[n+1]}\begin{bmatrix}2n\\n\end{bmatrix}$.

In this case the orthogonal polynomials whose moments are $C_n(q) = \frac{1}{[n+1]}\begin{bmatrix}2n\\n\end{bmatrix}$ have no known closed formulas. Therefore we consider instead the polynomials

$$f(n,x,q) = \sum_{j=0}^{\lfloor\frac{n}{2}\rfloor}(-1)^j q^{\binom{j}{2}}\begin{bmatrix}n-j\\j\end{bmatrix}x^{n-2j} \tag{82}$$

and the linear functional $\Lambda$ defined by $\Lambda\big(f(n,x,q)\big) = [n=0]$. As shown in [8] we also have

$$\Lambda\big(x^{2n}\big) = C_n(q) = \frac{1}{[n+1]}\begin{bmatrix}2n\\n\end{bmatrix} \text{ and } \Lambda\big(x^{2n+1}\big) = 0.$$

Therefore we get

$$\Lambda\big(f(2n,x,q)\big) = \Lambda\left(\sum_{j=0}^{n}(-1)^j q^{\binom{j}{2}}\begin{bmatrix}2n-j\\j\end{bmatrix}x^{2n-2j}\right) = \Lambda\left(\sum_{j=0}^{n}(-1)^{n-j} q^{\binom{n-j}{2}}\begin{bmatrix}n+j\\n-j\end{bmatrix}x^{2j}\right)$$

$$= \sum_{j=0}^{n}(-1)^{n-j} q^{\binom{n-j}{2}}\begin{bmatrix}n+j\\n-j\end{bmatrix}C_j(q) = [n=0],$$

which implies

$$\det\left(q^{\binom{i-j+1}{2}}\begin{bmatrix}i+j+1\\i-j+1\end{bmatrix}\right)_{i,j=0}^{n-1} = C_n(q) = \frac{1}{[n+1]}\begin{bmatrix}2n\\n\end{bmatrix}. \tag{83}$$

The first terms are

$$\det(1) = 1 = C_1(q), \quad \det\begin{pmatrix}1 & 1\\ q & [3]\end{pmatrix} = 1 + q^2 = C_2(q),$$

$$\det\begin{pmatrix}1 & 1 & 0\\ q & [3] & 1\\ q^3 & q+q^2+2q^3+q^4+q^5 & [5]\end{pmatrix} = (1-q+q^2)[5] = C_3(q),\cdots.$$

Since $\binom{i-j+1}{2} - \binom{i-j}{2} = i-j$

we also get

$$\det\left(q^{\binom{i-j}{2}}\begin{bmatrix}i+j+1\\i-j+1\end{bmatrix}\right)_{i,j=0}^{n-1}=C_n(q). \tag{84}$$

Let us also here consider the case $q=-1$. Since $C_n(q)=\begin{bmatrix}2n\\n\end{bmatrix}-q\begin{bmatrix}2n\\n-1\end{bmatrix}$ we get

$$C_{2n}(-1)=\begin{bmatrix}4n\\2n\end{bmatrix}_{-1}+\begin{bmatrix}4n\\2n-1\end{bmatrix}_{-1}=\binom{2n}{n}\text{ and }C_{2n+1}(-1)=\begin{bmatrix}4n+2\\2n+1\end{bmatrix}_{-1}+\begin{bmatrix}4n+2\\2n\end{bmatrix}_{-1}=\binom{2n+1}{n}$$

and thus $C_n(-1)=\left(\!\!\begin{array}{c}n\\ \left\lfloor\frac{n}{2}\right\rfloor\end{array}\!\!\right)$. This implies

$$\det\left((-1)^{\binom{i-j}{2}}\left(\begin{bmatrix}\frac{i+j+1}{2}\\j\end{bmatrix}\right)\right)_{i,j=0}^{n-1}=\left(\!\!\begin{array}{c}n\\ \left\lfloor\frac{n}{2}\right\rfloor\end{array}\!\!\right). \tag{85}$$

**Remark**

For $q=-1$ the polynomials $f(n,x,-1)=\sum_{j=0}^{\left\lfloor\frac{n}{2}\right\rfloor}(-1)^j(-1)^{\binom{j}{2}}\begin{bmatrix}n-j\\j\end{bmatrix}_{-1}x^{n-2j}$ are orthogonal with moments $\Lambda\left(x^{2n}\right)=\left(\!\!\begin{array}{c}n\\ \left\lfloor\frac{n}{2}\right\rfloor\end{array}\!\!\right)$.

As $q-$ analogue of (54) we get

**Theorem 7**

$$\det\left(q^{\binom{i-j}{2}}\begin{bmatrix}i+j+k\\i-j+1\end{bmatrix}\right)_{i,j=0}^{n-1}=\det\left(q^{\binom{i-j+1}{2}}\begin{bmatrix}i+j+k\\i-j+1\end{bmatrix}\right)_{i,j=0}^{n-1}=C_n^{(k)}(q)=\frac{[k]}{[2n+k]}\begin{bmatrix}2n+k\\n\end{bmatrix}.$$
(86)

**Proof**

It suffices to prove

$$\sum_{j=0}^{n}(-1)^{n-j}q^{\binom{n-j}{2}}\begin{bmatrix}n+j+k-1\\n-j\end{bmatrix}\frac{[k]}{[2j+k]}\begin{bmatrix}2j+k\\j\end{bmatrix}=[n=0]. \tag{87}$$

for $n,k\in\mathbb{N}$.

Using $(q;q)_{n-k} = \dfrac{(q;q)_n}{(q^{-n};q)_k}(-1)^k q^{\binom{k}{2}-nk}$ and $q-$Vandermonde's formula (cf. [12], II(6))

$$_2\varphi_1\left[\begin{matrix} a, q^{-n} \\ c \end{matrix}; q, q\right] = \sum_{k=0}^{n} \frac{(a;q)_k (q^{-n};q)_k}{(q;q)_k (c;q)_k} q^k = \frac{(c/a;q)_n}{(c;q)_n} a^n$$

we get

$$\sum_{j=0}^{n}(-1)^j q^{\binom{j}{2}} \begin{bmatrix} 2n-j+k-1 \\ j \end{bmatrix} \frac{[k]}{[2n-2j+k]} \begin{bmatrix} 2n-2j+k \\ n-j \end{bmatrix}$$

$$= \sum_{j=0}^{n}(-1)^j q^{\binom{j}{2}} \frac{(q;q)_{2n+k-1-j} [k] (q;q)_{2n-2j+k-1}}{(q;q)_j (q;q)_{2n+k-1-2j} (q;q)_{n-j} (q;q)_{n+k-j}}$$

$$= \sum_{j=0}^{n}(-1)^j q^{\binom{j}{2}} \frac{(q;q)_{2n+k-1-j} [k]}{(q;q)_j (q;q)_{n-j} (q;q)_{n+k-j}}$$

$$= \sum_{j=0}^{n}(-1)^j q^{\binom{j}{2}} \frac{(q;q)_{2n+k-1-j} [k]}{(q;q)_j (q;q)_{n-j} (q;q)_{n+k-j}}$$

$$= \frac{[k]}{[2n+k]}\begin{bmatrix} 2n+k \\ n \end{bmatrix} \sum_{j=0}^{n} q^j \frac{(q^{-n};q)_j (q^{-n-k};q)_j}{(q;q)_j (q^{-2n-k+1};q)_j} = \frac{[k]}{[2n+k]}\begin{bmatrix} 2n+k \\ n \end{bmatrix} {}_2\varphi_1\left[\begin{matrix} q^{-n}, q^{-n-k} \\ q^{-2n-k+1} \end{matrix}; q, q\right] = [n=0].$$

We could also imitate the proof of 2.1.3.2. Using $q-$Vandermonde's formula in the form

$$\sum_{k=0}^{r}\begin{bmatrix} r \\ k \end{bmatrix}\begin{bmatrix} s \\ n-k \end{bmatrix} q^{(r-k)(n-k)} = \begin{bmatrix} r+s \\ n \end{bmatrix} \quad \text{and} \quad \begin{bmatrix} -a \\ k \end{bmatrix} = \begin{bmatrix} a+k-1 \\ k \end{bmatrix}\left(-\frac{1}{q^a}\right)^k q^{-\binom{k}{2}}\ \text{we get}$$

$$\sum_{j=0}^{n}(-1)^{n-j} q^{\binom{n-j}{2}} \begin{bmatrix} n+j+k-1 \\ n-j \end{bmatrix} \frac{[k]}{[2j+k]} \begin{bmatrix} 2j+k \\ j \end{bmatrix} = (-1)^n \frac{[k]}{[n+k]} \sum_{j=0}^{n}(-1)^j q^{\binom{n-j}{2}} \begin{bmatrix} n+j+k-1 \\ n-j \end{bmatrix}\begin{bmatrix} n+k \\ n-j \end{bmatrix}$$

$$= (-1)^n \frac{[k]}{[n+k]} \sum_{j=0}^{n} \begin{bmatrix} -n-k \\ j \end{bmatrix}\begin{bmatrix} n+k \\ n-j \end{bmatrix} q^{j(j+k)} = \begin{bmatrix} n \\ 0 \end{bmatrix} = [n=0].$$

**Remark**

Formula (87) implies

$$\sum_{j=0}^{n}(-1)^{n-j} q^{\binom{n-j}{2}} \begin{bmatrix} (n+k)+(j+k) \\ (n+k)-(j+k) \end{bmatrix} C_j^{(2k+1)}(q)$$

$$= \sum_{j=0}^{n}(-1)^{n-j} q^{\binom{n-j}{2}} \begin{bmatrix} n+j+2k \\ n-j \end{bmatrix} \frac{[2k+1]}{[2j+2k+1]} \begin{bmatrix} 2j+2k+1 \\ j \end{bmatrix} = [n=0]$$

or $\sum_{j=0}^{n}(-1)^{n-j}q^{\binom{n-j}{2}}\begin{bmatrix}n+j\\n-j\end{bmatrix}C_{j-k}^{(2k+1)}(q)=[n=k].$

This implies the followin $q-$ analogue of (34):

$$\left((-q)^{\binom{i-j}{2}}\begin{bmatrix}i+j\\i-j\end{bmatrix}\right)^{-1}=\left(C_{i-j}^{(2j+1)}(q)\right). \tag{88}$$

As $q-$ analog of (55) we get

**Theorem 8**

Let $c(n,j)=q^{2\binom{n-j}{2}}\dfrac{\left(-q^{j};q\right)_{n-j}}{\left(-q;q\right)_{n-j}}\begin{bmatrix}j\\n-j\end{bmatrix}$. Then we get

$$\det\left(c(i+k+1,j+k)\right)_{i,j=0}^{n-1}=C(n,k,q)=\frac{[k]}{[2n+k]}\begin{bmatrix}2n+k\\n\end{bmatrix}\frac{(-q^{n+1};q)_{k-1}}{(-q;q)_{k-1}}. \tag{89}$$

It suffices to prove

**Lemma 9**

$$\sum_{j=0}^{n}(-1)^{n-j}c(n+k,j+k)\frac{[k]}{[2j+k]}\begin{bmatrix}2j+k\\j\end{bmatrix}\frac{(-q^{j+1};q)_{k-1}}{(-q;q)_{k-1}}=[n=0]. \tag{90}$$

For $k=1$ this has been found by George Andrews [2] and later a combinatorial proof has been given in [14].

To prove Lemma 9 we use the $q-$ Pfaff-Saalschütz summation (cf. [12],II.12)

$$_3\varphi_2\begin{bmatrix}a,b,q^{-n}\\c,abc^{-1}q^{1-n}\end{bmatrix};q,q\end{bmatrix}=\sum_{k=0}^{n}\frac{(a;q)_k(b;q)_k(q^{-n};q)_k}{(q;q)_k(c;q)_k(abc^{-1}q^{1-n};q)_k}q^k=\frac{(c/a;q)_n(c/b;q)_n}{(c;q)_n(c/ab;q)_n}$$

and the following elementary identities:

$$(-q^{n-k+1};q)_k=(q^{-n};q)_k q^{nk-\binom{k}{2}},\ (q;q)_{n-k}=\frac{(q;q)_n}{(q^{-n};q)_k}(-1)^k q^{\binom{k}{2}-nk},$$

$$(a;q)_k(-a;q)_k=(a^2;q^2)_k,\ (a;q)_{2k}=(a;q^2)_k(aq;q^2)_k.$$

First we change the order of summation

$$\sum_{j=0}^{n}(-1)^{n-j}q^{2\binom{n-j}{2}}\frac{(-q^{j+k};q)_{n-j}}{(-q;q)_{n-j}}\begin{bmatrix}j+k\\n-j\end{bmatrix}\frac{[k]}{[2j+k]}\begin{bmatrix}2j+k\\j\end{bmatrix}\frac{(-q^{j+1};q)_{k-1}}{(-q;q)_{k-1}}$$
$$=\sum_{j=0}^{n}(-1)^{j}q^{2\binom{j}{2}}\frac{(-q^{n-j+k};q)_{j}}{(-q;q)_{j}}\begin{bmatrix}n-j+k\\j\end{bmatrix}\frac{[k]}{[2n-2j+k]}\begin{bmatrix}2n-2j+k\\n-j\end{bmatrix}\frac{(-q^{n-j+1};q)_{k-1}}{(-q;q)_{k-1}}.$$

Then we make the following changes:

$$\left(-q^{n-j+k};q\right)_{j}=\left(-q^{-n-k+1};q\right)_{j}q^{(n+k-1)j-\binom{j}{2}}, \quad \frac{1}{(q;q)_{j}(q;q)_{n+k-2j}}=\frac{\left(q^{-n-k};q\right)_{2j}}{(q;q)_{j}(q;q)_{n+k}q^{\binom{2j}{2}-2j(n+k)}},$$

$$(q;q)_{2n+k-1-2j}=\frac{(q;q)_{2n+k-1}}{\left(q^{-2n-k+1};q\right)_{2j}}q^{\binom{2j}{2}-2j(2n+k-1)}, \quad \frac{1}{(q;q)_{n-j}}=\frac{\left(q^{-n};q\right)_{j}}{(q;q)_{n}}q^{nj-\binom{j}{2}},$$

$$\frac{1}{(q;q)_{n-j+k}}=\frac{\left(q^{-n-k};q\right)_{j}}{(q;q)_{n+k}}q^{(n+k)j-\binom{j}{2}}.$$

Observing that

$$\frac{\left(-q^{-n-k+1};q\right)_{j}\left(q^{-n};q\right)_{j}(-q^{n-j+1};q)_{k-1}}{(-q^{n+1};q)_{k-1}}=\frac{\left(q^{-2n};q^{2}\right)_{j}}{q^{(k-1)j}}$$

we get

$$\sum_{j=0}^{n}(-1)^{j}q^{2\binom{j}{2}}\frac{(-q^{n-j+k};q)_{j}}{(-q;q)_{j}}\begin{bmatrix}n-j+k\\j\end{bmatrix}\frac{[k]}{[2n-2j+k]}\begin{bmatrix}2n-2j+k\\n-j\end{bmatrix}\frac{(-q^{n-j+1};q)_{k-1}}{(-q;q)_{k-1}}$$
$$=\frac{[k]}{[2n+k]}\begin{bmatrix}2n+k\\n\end{bmatrix}\frac{(-q^{n+1};q)_{k-1}}{(-q;q)_{k-1}}q^{(k+1)j}\frac{\left(-q^{-n-k+1};q\right)_{j}\left(q^{-n-k};q\right)_{2j}\left(q^{-n};q\right)_{j}}{(-q;q)_{j}(q;q)_{j}\left(q^{-2n-k+1};q\right)_{2j}}\frac{(-q^{n-j+1};q)_{k-1}}{(-q^{n+1};q)_{k-1}}$$
$$=\frac{[k]}{[2n+k]}\begin{bmatrix}2n+k\\n\end{bmatrix}\frac{(-q^{n+1};q)_{k-1}}{(-q;q)_{k-1}}\frac{\left(q^{-2n};q^{2}\right)_{j}\left(q^{-n-k};q^{2}\right)_{j}\left(q^{-n-k+1};q^{2}\right)_{j}}{(q^{2};q^{2})_{j}\left(q^{-2n-k+1};q^{2}\right)_{j}\left(q^{-2n-k+2};q^{2}\right)_{j}}q^{2j}$$
$$=C(n,k,q)\,{}_3\varphi_2\begin{bmatrix}q^{-n-k},q^{-n-k+1},q^{-2n}\\q^{-2n-k+1},q^{-2n-k+2}\end{bmatrix};q^{2};q^{2}\end{bmatrix}=\frac{\left(q^{-n+1};q^{2}\right)_{n}\left(q^{-n};q^{2}\right)_{n}}{\left(q^{-2n-k+1};q^{2}\right)_{n}\left(q^{k};q^{2}\right)_{n}}=[n=0].$$

With Lemma 3 we get as $q-$analogue of (74)

**Theorem 10**

$$\det\left(q^{\binom{i-j+m}{2}}\begin{bmatrix}k+i+j+m\\i-j+m\end{bmatrix}\right)_{i,j=0}^{n-1} = q^{n\binom{m}{2}}\prod_{j=1}^{n}\frac{\prod_{\ell=0}^{2m-1}[2j-1+k+\ell]}{\prod_{\ell=0}^{m-1}[j+\ell][j+k+m+\ell]} = \det\left(C_{n-i+j}^{(2i+k+1)}(q)\right)_{i,j=0}^{m-1} \tag{91}$$

The proof is similar to the proof of (74).

A nice $q-$analogue of (45) is

$$\det\left(q^{\binom{i-j}{2}}\frac{[2i+k+1]}{[i+j+k]}\begin{bmatrix}i+j+k\\i-j+1\end{bmatrix}\right)_{i,j=0}^{n-1} = \begin{bmatrix}2n+k-1\\n\end{bmatrix}, \tag{92}$$

which follows from

$$\sum_{j=0}^{n}(-1)^j q^{\binom{j}{2}}\frac{[2n+k-1]}{[2n-j+k-1]}\begin{bmatrix}2n-j+k-1\\j\end{bmatrix}\begin{bmatrix}2n-2j+k-1\\n-j\end{bmatrix} =$$

$$\begin{bmatrix}2n+k-1\\n\end{bmatrix}{}_2\varphi_1\left[\begin{matrix}q^{-n},q^{-n-k+1}\\q^{-2n-k+2}\end{matrix};q,q\right] = \frac{(q^{-n+1};q)_n}{(q^{-2n-k+2};q)_n} = [n=0].$$

**Remark**

For $k=1$ this is also a consequence of the fact (cf. [8]) that the polynomials

$$L_n(x,q) = \sum_{j=0}^{\lfloor\frac{n}{2}\rfloor}(-1)^j q^{\binom{j}{2}}\frac{[n]}{[n-j]}\begin{bmatrix}n-j\\j\end{bmatrix}x^{n-2j} \text{ satisfy } \sum_{k=0}^{\lfloor\frac{n}{2}\rfloor}\begin{bmatrix}n\\k\end{bmatrix}L_{n-2k}(x,q) = x^n.$$

This implies that the linear functional $\Lambda$, defined by $\Lambda(L_n(x,q)) = [n=0]$ satisfies

$$\Lambda(x^{2n}) = \begin{bmatrix}2n\\n\end{bmatrix} \text{ and } \Lambda(x^{2n+1}) = 0.$$

Since

$$L_{2n}(x,q) = \sum_{j=0}^{n}(-1)^j q^{\binom{j}{2}}\frac{[2n]}{[2n-j]}\begin{bmatrix}2n-j\\j\end{bmatrix}x^{2n-2j} = \sum_{j=0}^{n}(-1)^{n-j}q^{\binom{n-j}{2}}\frac{[2n]}{[n+j]}\begin{bmatrix}n+j\\n-j\end{bmatrix}x^{2j}$$

we get by applying $\Lambda$

$$\sum_{j=0}^{n}(-1)^{n-j}q^{\binom{n-j}{2}}\frac{[2n]}{[n+j]}\begin{bmatrix}n+j\\n-j\end{bmatrix}\begin{bmatrix}2j\\j\end{bmatrix}=[n=0]$$

By Lemma 1 this implies (92).

Let us also consider the polynomials $\ell_n(x,q)=\sum_{j=0}^{\lfloor\frac{n}{2}\rfloor}(-1)^j q^{2\binom{j}{2}}\frac{[n]}{[n-j]}\begin{bmatrix}n-j\\j\end{bmatrix}x^{n-2j}$ and the linear functional $\lambda$ defined by $\lambda\bigl(\ell_n(x,q)\bigr)=[n=0]$. The corresponding moments $\lambda\bigl(x^n\bigr)$ are
$1,\ 0,\ 1+q,\ 0,\ \bigl(1+q^2\bigr)\bigl(1+2q\bigr),\ 0,\ \bigl(1+q^3\bigr)\bigl(1+3q+3q^2+3q^3\bigr),\ 0,$
$\bigl(1+q^4\bigr)\bigl(1+2q+2q^3\bigr)\bigl(1+2q+2q^2+2q^3\bigr),\cdots.$

As far as I know there is no closed formula for these moments. But it seems that all coefficients are positive. More generally it seems that

$$\det\left(q^{2\binom{i-j}{2}}\frac{[2i+k+1]}{[i+j+k]}\begin{bmatrix}i+j+k\\i-j+1\end{bmatrix}\right)_{i,j=0}^{n-1}=b(n,k,q),\tag{93}$$

where all non-vanishing coefficients of $b(n,j,q)$ are positive.

It would be interesting if there is also a nice $q-$analogue of (46).

A companion to Theorem 10 is

**Theorem 11**

*Let*

$$B_n(x,m,q)=\left(q^{\binom{i-j+m}{2}}\frac{[2i+x+2m-1]}{[i+j+x+m-1]}\begin{bmatrix}i+j+x+m-1\\i-j+m\end{bmatrix}\right)_{i,j=0}^{n-1}\tag{94}$$

*and*

$$H_m(x,n,q)=\bigl(h_{n-i+j}(x+2i,q)\bigr)_{i,j=0}^{m-1}\tag{95}$$

*with* $h_n(x,q)=\begin{bmatrix}2n+x-1\\n\end{bmatrix}$. *Then*

$$\det B_n(x,m,q)=q^{n\binom{m}{2}}\frac{1}{(1-q)^{mn}}\prod_{j=0}^{m-1}\frac{[j]!}{[n+j]!}\prod_{j=1}^{n}\bigl(q^{x+2j-2};q\bigr)_{m-j}\bigl(q^{x+2m+j-2};q\bigr)_j=\det H_m(x,n,q).$$
(96)

**Proof**

Let $w(n,x,m) = q^{n\binom{m}{2}} \dfrac{1}{(1-q)^{mn}} \prod_{j=0}^{m-1} \dfrac{[j]!}{[n+j]!} \prod_{j=1}^{n} \left(q^{x+2j-2};q\right)_{m-j} \left(q^{x+2m+j-2};q\right)_{j}.$

We have $w(0,x,m) = w(n,x,0) = 1$ and $\det B_0(x,m,q) = \det H_0(x,n,q) = 1.$

Further we have

$$w(n,x,1) = \begin{bmatrix} 2n+x-1 \\ n \end{bmatrix} \tag{97}$$

and

$$w(1,x,m) = q^{\binom{m}{2}} \begin{bmatrix} x+m-1 \\ m \end{bmatrix} \frac{[x+2m-1]}{[x+m-1]}. \tag{98}$$

By (92) and (97) we have $\det B_n(x,1,q) = w(n,x,1).$

Identity (98) follows from

$$w(1,x,m) = q^{\binom{m}{2}} \frac{1}{(1-q)^m} \prod_{j=0}^{m-1} \frac{[j]!}{[1+j]!} \left(q^x;q\right)_{m-1} \left(q^{x+2m-1};q\right)_1$$
$$= q^{\binom{m}{2}} \frac{[x]\cdots[x+m-2][x+2m-1]}{[m]!} == q^{\binom{m}{2}} \begin{bmatrix} x+m-1 \\ m \end{bmatrix} \frac{[x+2m-1]}{[x+m-1]}.$$

Now we must show that

$$\det H_m(x,1,q) = \det\left(h_{1-i+j}(x+2i,q)\right)_{i,j=0}^{m-1} = w(1,x,m).$$

It remains to prove that

$$\det\left(\begin{bmatrix} 2i+x+1 \\ i-j+1 \end{bmatrix}\right)_{i,j=0}^{m-1} = q^{\binom{m}{2}} \begin{bmatrix} x+m-1 \\ m \end{bmatrix} \frac{[x+2m-1]}{[x+m-1]}.$$

By Lemma 1 it suffices to prove that

$$\sum_{j=0}^{n} (-1)^{n-j} \begin{bmatrix} 2n+x-1 \\ n-j \end{bmatrix} q^{\binom{j}{2}} \begin{bmatrix} x+j-1 \\ j \end{bmatrix} \frac{[x+2j-1]}{[x+j-1]} = [n=0].$$

This reduces to

$$\begin{bmatrix} 2n+x-1 \\ n \end{bmatrix} {}_4\varphi_3 \left[\begin{matrix} q^{x-1}, q^{\frac{x+1}{2}}, -q^{\frac{x+1}{2}}, q^{-n} \\ q^{\frac{x-1}{2}}, -q^{\frac{x-1}{2}}, q^{n+x} \end{matrix}; q, q^n\right] = [n=0],$$

which is the special case of identity [12], (II.1)

$$_4\varphi_3\left[\begin{array}{c} a,-q\sqrt{a},b,q^{-n} \\ -\sqrt{a},\dfrac{aq}{b},aq^{n+1} \end{array};q,\dfrac{q^{n+1}\sqrt{a}}{b}\right] = \dfrac{(aq;q)_n\left(\dfrac{q\sqrt{a}}{b};q\right)_n}{\left(q\sqrt{a};q\right)_n\left(\dfrac{aq}{b};q\right)_n}$$

for $a = q^{x-1}$, $b = q^{\frac{x+1}{2}}$.

To prove (96) for all values $m, n$ we use the condensation method. In our case this reduces to

$$\det B_n(x,m,q)\det B_{n-2}(x+2,m,q) - \det B_{n-1}(x+2,m,q)\det B_{n-1}(x,m,q)$$
$$+ \det B_{n-1}(x,m+1,q)\det B_{n-1}(x+2,m-1,q) = 0 \qquad (99)$$

and to

$$\det H_m(x,n,q)\det H_{m-2}(x+2,n,q) - \det H_{m-1}(x+2,n,q)\det H_{m-1}(x,n,q)$$
$$+ \det H_{m-1}(x,n+1,q)\det H_{m-1}(x+2,n-1,q) = 0. \qquad (100)$$

Since

$$\det B_n(x,m,q)\det B_{n-2}(x+2,m,q) : \det B_{n-1}(x+2,m,q)\det B_{n-1}(x,m,q)$$
$$: \det B_{n-1}(x,m+1,q)\det B_{n-1}(x+2,m-1,q)$$
$$= [n-1]\left(1-q^{n+2m-2+x}\right) : [n+m-1]\left(1-q^{n+m-2+x}\right) : q^{n-1}[m]\left(1-q^{x+m-1}\right)$$

and $[n-1]\left(1-q^{n+2m-2+x}\right) - [n+m-1]\left(1-q^{n+m-2+x}\right) + q^{n-1}[m]\left(1-q^{x+m-1}\right) = 0$

we get (99).

Since

$$\det H_m(x,n,q)\det H_{m-2}(x+2,n,q) : \det H_{m-1}(x+2,n,q)\det H_{m-1}(x,n,q)$$
$$: \det H_{m-1}(x,n+1,q)\det H_{m-1}(x+2,n-1,q)$$
$$= q^n[m-1]\left(1-q^{x+m-2)}\right) : [m+n-1]\left(1-q^{x+m+n-2)}\right) : [n]\left(1-q^{x+2m+n}\right)$$

and $q^n[m-1]\left(1-q^{x+m-2)}\right) - [m+n-1]\left(1-q^{x+m+n-2)}\right) + [n]\left(1-q^{x+2m+n}\right) = 0$

we get (100).

### 3.3. q-Chebyshev polynomials and the q-Catalan numbers of George Andrews

For the monic Chebyshev polynomials of the second kind

$$u_n(x) = \sum_{k=0}^{\lfloor n/2 \rfloor} \left(-\frac{1}{4}\right)^k \binom{n-k}{k} x^{n-2k} = \frac{F_n(2x)}{2^n}$$

there exist $q$-analogues with nice formulas of both the orthogonal polynomials and its moments.

As analogues of $u_n(x)$ we get the monic $q$-Chebyshev polynomials of the second kind (cf.[7])

$$u_n(x,q) = \sum_{k=0}^{\lfloor n/2 \rfloor} \frac{(-1)^k}{(-q;q)_k \left(-q^{n+1-k};q\right)_k} q^{k^2} \begin{bmatrix} n-k \\ k \end{bmatrix} x^{n-2k} \qquad (101)$$

which satisfy $u_n(x,q) = x u_{n-1}(x,q) - \dfrac{q^{n-1}}{\left(1+q^{n-1}\right)\left(1+q^n\right)} u_{n-2}(x,q)$ with initial values $u_0(x,q) = 1$ and $u_1(x,q) = x$. Their moments are Andrews' $q$-Catalan numbers

$$M_n = \frac{1}{[n+1]} \begin{bmatrix} 2n \\ n \end{bmatrix} \frac{1+q}{1+q^{n+1}} \frac{q^n}{(-q;q)_n^2}. \qquad (102)$$

Here we get as analogue of (54)

$$\det\left(\frac{q^{(i+1-j)^2}}{(-q;q)_{i+1-j}\left(-q^{i+j+k+1};q\right)_{i+1-j}} \begin{bmatrix} i+j+k \\ i+1-j \end{bmatrix}\right)_{i,j=0}^{n-1} = q^n \frac{1+q^k}{1+q^{n+k}} \frac{[k]}{[2n+k]} \begin{bmatrix} 2n+k \\ n \end{bmatrix} \frac{1}{(-q;q)_n \left(-q^k;q\right)_n}.$$

The proofs use orthogonality and will be omitted.

### 4. Concluding remarks

In [10] we studied the sequence $a(n) = C_n \mod 2$ as a sequence of numbers from $\{0,1\}$ and obtained their orthogonal polynomials $p_n(x)$. Let $p_{2n}\left(\sqrt{x}\right) = \sum_{j=0}^{n} (-1)^{n-j} p(n,j) x^j$. The first terms of the sequence $p_{2n}\left(\sqrt{x}\right)$ are $1,\ x-1,\ x^2+x-1,\ x^3+x^2-1,\ x^4+x^3+x^2-1,\cdots$, and the first terms of the matrix $\left(p(i,j)\right)_{i,j\geq 0}$ are

$$\begin{pmatrix} 1 & & & & & \\ 1 & 1 & & & & \\ -1 & -1 & 1 & & & \\ 1 & 0 & -1 & 1 & & \\ 1 & 0 & 1 & -1 & 1 & \\ -1 & -1 & 1 & 0 & 1 & 1 \end{pmatrix}$$

Since by [10], (0.1), $\det\left(a(i+j)\right)_{i,j=0}^{n-1} = (-1)^{\binom{n}{2}}$ we get from (68)

$$\det\left(p(i+m,j)\right)_{i,j=0}^{n-1} = (-1)^{\binom{m}{2}} \det\left(a(i+j+n)\right)_{i,j=0}^{m-1}. \tag{103}$$

As special case we get

$$\det\left(p(i+1,j)\right)_{i,j=0}^{n-1} = C_n \bmod 2. \tag{104}$$

Consider now the polynomials

$$r_n(x) = \sum_{j=0}^{n} (-1)^{n-j} \binom{n+j}{n-j} (\bmod 2) x^j = \sum_{j=0}^{n} (-1)^{n-j} r(n,j) x^j \tag{105}$$

which also satisfy $r_n(x) \equiv p_{2n}\left(\sqrt{x}\right) \bmod 2$.

The first terms are $1,\ x-1,\ x^2-x+1,\ x^3-x^2-1,\ x^4-x^3+x^2+1,\cdots$.

Here we have

$$\sum_{j=0}^{n} (-1)^{n-j} r(n,j) \left(C_j \bmod 2\right) = [n=0]. \tag{106}$$

Therefore we also get

$$\det\left(r(i+1,j)\right)_{i,j=0}^{n-1} = \det\left(\binom{i+j+1}{i-j+1} \bmod 2\right)_{i,j=0}^{n-1} = C_n \bmod 2. \tag{107}$$

For the proof of (106) observe that for $n > 0$ we have

$$\sum_{j=0}^{n} (-1)^j r(n,j)\left(C_j \bmod 2\right) = \sum_{j=0}^{n} (-1)^{2^j-1} r(n, 2^j-1) = r(n,0) - \sum_{j=1}^{n} r(n, 2^j-1)$$

$$= 1 - \sum_{j=1}^{n} \binom{n-1+2^j}{n+1-2^j} \bmod 2 = 0,$$

because there is exactly one $j > 0$ such that $\binom{n-1+2^j}{n+1-2^j} \equiv 1 \bmod 2$.

Since by Lucas's theorem $\binom{2a}{2b} \equiv \binom{2a+1}{2b+1} \equiv \binom{a}{b} \bmod 2$ it suffices to show that for each non-negative integer $m$ there is exactly one $j \geq 0$ such that $\binom{m+2^j}{m+1-2^j} \equiv 1 \bmod 2$.

To prove this last fact we use Lucas' theorem: Assume $a = \sum_k a_k 2^k$, $b = \sum_k b_k 2^k$ are the representations of $a, b$ in base 2. Then $\binom{a}{b} \equiv \prod_k \binom{a_k}{b_k} \bmod 2$.

Let $m = m_0 + 2m_1 + 2^2 m_2 + \cdots$ be the representation of $m$ in base 2 and let $h$ be the smallest $j > 0$ such that $m_j = 0$. Then

$$m + 2^h = 1 + 2 + \cdots + 2^h + 2^{h+1} m_{h+1} + \cdots \text{ and } m + 1 - 2^h = 2^{h+1} m_{h+1} + \cdots.$$

By Lucas's theorem we get $\binom{m+2^h}{m+1-2^h} \equiv 1 \bmod 2$, since all $\binom{a_k}{b_k}$ are $\binom{1}{0} = 1$ or $\binom{m_j}{m_j} = 1$.

For $j < h$ we get $m + 2^j = 1 + \cdots + 2^{j-1} + 0 \cdot 2^j + \cdots$ and $m + 1 - 2^j = 2^j + \cdots$. Therefore we have $\binom{a_j}{b_j} = \binom{0}{1} = 0$.

For $j > h$ and $m + 1 - 2^j \geq 0$ we have $(m+2^j)_h = 0$ and $(m+1-2^j)_h = 1$ and therefore also the factor $\binom{0}{1} = 0$.

It seems that more than Formula (107) is true:

**Conjecture 12**

$$\left( (-1)^{i-j} \binom{i+j}{i-j} \bmod 2 \right)^{-1} = \left( C_{i-j}^{(2j+1)} \bmod 2 \right)_{i,j \geq 0}. \tag{108}$$

## Conjecture 13

*For $n, m, k > 0$ the following identities hold:*

$$\det\left(\binom{i+j+k}{i-j+1} \mod 2\right)_{i,j=0}^{n-1} = (-1)^{k-1} C_n^{(k)} \mod 2, \qquad (109)$$

$$\det\left(\binom{i+j+m}{i-j+m} \mod 2\right)_{i,j=0}^{n-1} = \det\left(C_{n-i+j}^{(2i+1)} \mod 2\right)_{i,j=0}^{m-1}. \qquad (110)$$

For other moduli I have only found one interesting fact. Let $\mu(n)$ denote the residue modulo 3 where the residues belong to $\{-1, 0, 1\}$, i.e. $\mu(3n) = 0$, $\mu(3n+1) = 1$ and $\mu(3n+2) = -1$. Then we get

## Conjecture 14

$$\det\left(\binom{i+j+1}{i-j+1} \mod 3\right)_{i,j=0}^{n-1} = \mu(C_n). \qquad (111)$$

Let me also mention $q$-analogues of (36) and (47)

## Theorem 15

Let $A_n(x,q) = \left(q^{\binom{i-j}{2}} \begin{bmatrix} i+j+x \\ i-j+1 \end{bmatrix}\right)_{i,j=0}^{n-1}$ and $n+1 \le m \le 2n-1$. Then $A_n(-m, q) v_{n,m} = 0$ if

$$v_{n,m} = \left( q^{\frac{\left(\lfloor\frac{m}{2}\rfloor - j\right)\left(\lfloor\frac{m}{2}\rfloor + 5 - 3j\right)}{2}} \frac{[m]}{[m-j]} \begin{bmatrix} m-j \\ j \end{bmatrix} \right)^t_{0 \le j \le n-1} \quad \text{for even } m \text{ and}$$

$$v_{n,m} = \left( q^{\frac{\left(\lfloor\frac{m}{2}\rfloor - j\right)\left(\lfloor\frac{m}{2}\rfloor + 7 - 3j\right)}{2}} \frac{[m]}{[m-j]} \begin{bmatrix} m-j \\ j \end{bmatrix} \right)^t_{0 \le j \le n-1} \quad \text{for odd } m.$$

If $B_n(x,q) = \left(q^{\binom{i-j}{2}} \frac{[2i+x+1]}{[i+j+x]} \begin{bmatrix} i+j+x \\ i-j+1 \end{bmatrix}\right)_{i,j=0}^{n-1}$ and $n \le m \le 2n-1$. Then

$B_n(-m, q) u_{n,m} = 0$ if

$$u_{n,m} = \left( q^{\frac{\left[\left\lfloor\frac{m}{2}\right\rfloor - j\right]\left[\left\lfloor\frac{m}{2}\right\rfloor + 3 - 3j\right]}{2}} \begin{bmatrix} m-j \\ j \end{bmatrix} \right)^t_{0 \le j \le n-1} \quad \text{for even } m \text{ and}$$

$$u_{n,m} = \left( q^{\frac{\left[\left\lfloor\frac{m}{2}\right\rfloor - j\right]\left[\left\lfloor\frac{m}{2}\right\rfloor + 5 - 3j\right]}{2}} \begin{bmatrix} m-j \\ j \end{bmatrix} \right)^t_{0 \le j \le n-1} \quad \text{for odd } m.$$

This follows from

**Lemma 16**

*Let $x$ be an indeterminate. Then the following identities hold:*

$$\sum_{j=0}^{i+1} q^{\binom{i-j}{2}} \begin{bmatrix} i+j-2x \\ i-j+1 \end{bmatrix} q^{\frac{(x-j)(x+5-3j)}{2}} \frac{1-q^{2x}}{1-q^{2x-j}} \begin{bmatrix} 2x-j \\ j \end{bmatrix} = 0, \tag{112}$$

$$\sum_{j=0}^{i+1} q^{\binom{i-j}{2}} \begin{bmatrix} i+j-2x-1 \\ i-j+1 \end{bmatrix} q^{\frac{(x-j)(x+7-3j)}{2}} \frac{1-q^{2x+1}}{1-q^{2x+1-j}} \begin{bmatrix} 2x+1-j \\ j \end{bmatrix} = 0, \tag{113}$$

$$\sum_{j=0}^{i+1} q^{\binom{i-j}{2}} \frac{1}{[i+j-2x]} \begin{bmatrix} i+j-2x \\ i-j+1 \end{bmatrix} q^{\frac{(x-j)(x+3-3j)}{2}} \begin{bmatrix} 2x-j \\ j \end{bmatrix} = 0, \tag{114}$$

$$\sum_{j=0}^{i+1} q^{\binom{i-j}{2}} \frac{1}{[i+j-2x-1]} \begin{bmatrix} i+j-2x-1 \\ i-j+1 \end{bmatrix} q^{\frac{(x-j)(x+5-3j)}{2}} \begin{bmatrix} 2x+1-j \\ j \end{bmatrix} = 0. \tag{115}$$

Let us deduce the first assertion of Theorem 14. The other assertions follow in a similar way. We must show that for $n+1 \le 2m \le 2n-1$

$$\sum_{j=0}^{n-1} q^{\binom{i-j}{2}} \begin{bmatrix} i+j-2m \\ i-j+1 \end{bmatrix} q^{\frac{(m-j)(m+5-3j)}{2}} \frac{1-q^{2m}}{1-q^{2m-j}} \begin{bmatrix} 2m-j \\ j \end{bmatrix} = 0$$

for $0 \le i \le n-1$.

By (112) this is true for $i < n-1$.

For $i = n-1$ we also have by (112)

$$\sum_{j=0}^{n} q^{\binom{n-1-j}{2}} \begin{bmatrix} n-1+j-2x \\ n-j \end{bmatrix} q^{\frac{(x-j)(x+5-3j)}{2}} \frac{1-q^{2x}}{1-q^{2x-j}} \begin{bmatrix} 2x-j \\ j \end{bmatrix} = 0$$

and therefore

$$\sum_{j=0}^{n-1} q^{\binom{n-1-j}{2}} \begin{bmatrix} n-1+j-2x \\ n-j \end{bmatrix} q^{\frac{(x-j)(x+5-3j)}{2}} \frac{1-q^{2x}}{1-q^{2x-j}} \begin{bmatrix} 2x-j \\ j \end{bmatrix} = -q^{\frac{(x-n)(x+5-3n)}{2}} \frac{1-q^{2x}}{1-q^{2x-n}} \begin{bmatrix} 2x-n \\ n \end{bmatrix}.$$

This vanishes for $n+1 \leq 2x \leq 2n-1$.

**Proof of Lemma 16**

Since the proofs of all assertions are almost identical we give only the proof of (112) in detail.

Using the identities $\begin{bmatrix} 2x-j \\ j \end{bmatrix} = \begin{bmatrix} 2j-2x-1 \\ j \end{bmatrix} (-1)^j q^{j(2x-j)} q^{-\binom{j}{2}}$,

$$\frac{1}{(q;q)_{i-j+1}} = (-1)^j q^{ij - \frac{j^2-3j}{2}} \frac{(q^{-i-1};q)_j}{(q;q)_{i+1}},$$

$$\begin{bmatrix} i+j-2x \\ i-j+1 \end{bmatrix} \begin{bmatrix} 2j-2x-1 \\ j \end{bmatrix} \frac{1-q^{2x}}{1-q^{2x-j}} = q^j \frac{1-q^{-2x}}{1-q^{j-2x}} \frac{[i+j-2x]\cdots[2j-2x][2j-2x-1]\cdots[j-2x]}{[i-j+1]![j]!}$$

$$= q^j \frac{(1-q^{-2x})\cdots(1-q^{-2x+i})(1-q^{-2x+i+1})\cdots(1-q^{-2x+i+j})}{(1-q^{-2x+1})\cdots(1-q^{-2x+j})} \frac{1}{(q;q)_j (q;q)_{i-j+1}}$$

$$= q^j \frac{(q^{1+i-2x};q)_j (q^{-2x};q)_{1+i}}{(q^{1-2x};q)_j} \frac{1}{(q;q)_j (q;q)_{i-j+1}},$$

and the $q$ – Chu-Vandermonde formula ${}_2\varphi_1 \begin{bmatrix} q^{-n}, a \\ c \end{bmatrix}; q, q \end{bmatrix} = \frac{\left(\frac{c}{a};q\right)_n}{(c;q)_n} a^n$ (cf. [12],II.6), we get

$$\sum_{j=0}^{i+1} q^{\binom{i-j}{2}} \begin{bmatrix} i+j-2x \\ i-j+1 \end{bmatrix} q^{\frac{(x-j)(x+5-3j)}{2}} \frac{1-q^{2x}}{1-q^{2x-j}} \begin{bmatrix} 2x-j \\ j \end{bmatrix}$$

$$= q^{\binom{i}{2} + \frac{x^2+5x}{2}} \sum_{j=0}^{i+1} q^{-ij-j+2j^2-2jx} \begin{bmatrix} i+j-2x \\ i-j+1 \end{bmatrix} \frac{1-q^{2x}}{1-q^{2x-j}} \begin{bmatrix} 2j-2x-1 \\ j \end{bmatrix} (-1)^j q^{j(2x-j)} q^{-\binom{j}{2}}$$

$$= q^{\binom{i}{2} + \frac{x^2+5x}{2}} \sum_{j=0}^{i+1} (-1)^j q^{-ij - \frac{j}{2} + \frac{j^2}{2}} \frac{(q^{1+i-2x};q)_j (q^{-2x};q)_{1+i}}{(q^{1-2x};q)_j} \frac{1}{(q;q)_j (q;q)_{i-j+1}}$$

$$= q^{\binom{i}{2} + \frac{x^2+5x}{2}} \frac{(q^{-2x};q)_{1+i}}{(q;q)_{i+1}} \sum_{j=0}^{i+1} q^j \frac{(q^{-i-1};q)_j (q^{1+i-2x};q)_j (q^{-2x};q)_{1+i}}{(q^{1-2x};q)_j (q;q)_j}$$

$$= q^{\binom{i}{2}+\frac{x^2+5x}{2}} \frac{\left(q^{-2x};q\right)_{1+i}}{(q;q)_{i+1}} {}_2\varphi_1\left[\begin{matrix} q^{-i-1}, q^{1+i-2x} \\ q^{1-2x} \end{matrix}; q, q\right] = q^{\binom{i}{2}+\frac{x^2+5x}{2}} \frac{\left(q^{-2x};q\right)_{1+i}}{(q;q)_{i+1}} \frac{\left(q^{-i};q\right)_{i+1}}{\left(q^{1-2x};q\right)_{i+1}} q^{(1+i-2x)(1+i)}$$

$$= q^{1+\frac{3i}{2}+\frac{3i^2}{2}+\frac{x}{2}-2ix+\frac{x^2}{2}} \frac{\left(q^{-2x};q\right)_{1+i}\left(q^{-i};q\right)_{i+1}}{(q;q)_{i+1}\left(q^{1-2x};q\right)_{i+1}} = 0.$$

In the same manner we obtain

$$\sum_{j=0}^{i+1} q^{\binom{i-j}{2}} \begin{bmatrix} i+j-2x-1 \\ i-j+1 \end{bmatrix} q^{\frac{(x-j)(x+7-3j)}{2}} \frac{1-q^{2x+1}}{1-q^{2x+1-j}} \begin{bmatrix} 2x+1-j \\ j \end{bmatrix} = q^{\binom{i}{2}+\frac{x^2+7x}{2}} \frac{\left(q^{-2x-1};q\right)_{1+i}}{(q;q)_{i+1}} {}_2\varphi_1\left[\begin{matrix} q^{-i-1}, q^{i-2x} \\ q^{-2x} \end{matrix}; q, q\right]$$

$$= q^{\frac{i}{2}+\frac{3i^2}{2}+\frac{3x}{2}-2ix+\frac{x^2}{2}} \frac{\left(q^{-1-2x};q\right)_{1+i}\left(q^{-i};q\right)_{i+1}}{(q;q)_{i+1}\left(q^{-2x};q\right)_{i+1}} = 0,$$

$$\sum_{j=0}^{i+1} q^{\binom{i-j}{2}} \frac{1}{[i+j-2x]} \begin{bmatrix} i+j-2x \\ i-j+1 \end{bmatrix} q^{\frac{(x-j)(x+3-3j)}{2}} \begin{bmatrix} 2x-j \\ j \end{bmatrix} = q^{\binom{i}{2}+\frac{x^2+3x}{2}} \frac{\left(q^{-2x};q\right)_i}{(q^2;q)_i} {}_2\varphi_1\left[\begin{matrix} q^{-i-1}, q^{i-2x} \\ q^{-2x} \end{matrix}; q, q\right]$$

$$= \frac{q^{\frac{i}{2}+\frac{3i^2}{2}+\frac{x^2}{2}-\frac{x}{2}-2xi}\left(q^{-2x};q\right)_i\left(q^{-i};q\right)_{i+1}}{(q^2;q)_i\left(q^{-2x};q\right)_{i+1}} = 0$$

and

$$\sum_{j=0}^{i+1} q^{\binom{i-j}{2}} \frac{1}{[i+j-2x-1]} \begin{bmatrix} i+j-2x-1 \\ i-j+1 \end{bmatrix} q^{\frac{(x-j)(x+5-3j)}{2}} \begin{bmatrix} 2x+1-j \\ j \end{bmatrix}$$

$$= q^{\binom{i}{2}+\frac{x^2+5x}{2}} \frac{\left(q^{-2x-1};q\right)_i}{(q^2;q)_i} {}_2\varphi_1\left[\begin{matrix} q^{-i-1}, q^{i-2x-1} \\ q^{-2x-1} \end{matrix}; q, q\right] = \frac{q^{1-\frac{i}{2}+\frac{3i^2}{2}+\frac{x^2}{2}+\frac{x}{2}-2xi}\left(q^{-2x-1};q\right)_i\left(q^{-i};q\right)_{i+1}}{(q^2;q)_i\left(q^{-2x-1};q\right)_{i+1}} = 0.$$

**Remark**

Theorem 10 and Theorem 11 are of the form

$$\det\left(p(i+m,j,x)\right)_{i,j=0}^{n-1} = \det\left(h_{n-i+j}(x+r(i,j,m))\right)_{i,j=0}^{m-1} \text{ for some } r(i,j,m).$$

There are also other examples such as

$$\det\left(q^{\binom{i-j+m}{2}}\begin{bmatrix} i+x+m \\ i-j+m \end{bmatrix}\right)_{i,j=0}^{n-1} = \det\left(h_{n-i+j}(x+m-1)\right)_{i,j=0}^{m-1} = q^{n\binom{m}{2}} \prod_{j=0}^{n-1} \frac{\left(q^{x+j+1};q\right)_{m+n-1-2j}}{\left(q^{j+1};q\right)_{m+n-1-2j}},$$

where $h_n(x) = \begin{bmatrix} n+x \\ n \end{bmatrix}$.

It would be interesting if there is a similar characterization of such pairs of $n \times n-$ determinants and corresponding $m \times m -$ determinants as in Theorem 5.